\documentclass[final,1p,times]{elsarticle}

\usepackage{amsmath}
\usepackage{amssymb}
\usepackage{amsthm}
\allowdisplaybreaks

\usepackage{hyperref}

\usepackage{mathtools}

\usepackage{mismath}

\usepackage{bm} 
\usepackage{bbm}

\usepackage{graphicx}
\graphicspath{{figs/}} 

\usepackage{tikz}
\usetikzlibrary{math}
\usetikzlibrary{calc}

\usepackage{subcaption}

\usepackage{multirow}

\usepackage{siunitx}
\usepackage{makecell}
\setcellgapes{3pt}

\usepackage{algorithm}

\usepackage{algpseudocode}

\usepackage[capitalize]{cleveref}

\usepackage{easyReview}
\setreviewsoff

\newcommand{\TD}{\Gamma_{\rm{D}}}
\newcommand{\TN}{\Gamma_{\rm{N}}}

\newcommand{\TC}{\Gamma_{\rm{C}}}

\DeclareMathOperator{\Div}{div}

\DeclareMathOperator{\Span}{span}

\DeclarePairedDelimiter{\CurlyBrackets}{\{}{\}}

\newtheorem{theorem}{Theorem}[section]
\newtheorem{lemma}[theorem]{Lemma}

\theoremstyle{definition}
\newtheorem{definition}{Definition}[section]
\newtheorem{assumption}[theorem]{Assumption}

\theoremstyle{remark}

\definecolor{MyColor1}{HTML}{CECE5A}
\definecolor{MyColor2}{HTML}{C51605}
\definecolor{MyColor3}{HTML}{B7B7B7}
\definecolor{MyColor4}{HTML}{FFE17B}

\usepackage{array}
\newcolumntype{L}[1]{>{\raggedright\let\newline\\\arraybackslash\hspace{0pt}}m{#1}}
\newcolumntype{C}[1]{>{\centering\let\newline\\\arraybackslash\hspace{0pt}}m{#1}}
\newcolumntype{R}[1]{>{\raggedleft\let\newline\\\arraybackslash\hspace{0pt}}m{#1}}

\definecolor{MyColor1}{HTML}{CECE5A}
\definecolor{MyColor2}{HTML}{C51605}
\definecolor{MyColor3}{HTML}{B7B7B7}
\definecolor{MyColor4}{HTML}{FFE17B}

\usepackage{array}
\newcolumntype{L}[1]{>{\raggedright\let\newline\\\arraybackslash\hspace{0pt}}m{#1}}
\newcolumntype{C}[1]{>{\centering\let\newline\\\arraybackslash\hspace{0pt}}m{#1}}
\newcolumntype{R}[1]{>{\raggedleft\let\newline\\\arraybackslash\hspace{0pt}}m{#1}}

\journal{arXiv}

\begin{document}

\begin{frontmatter}

  \title{An iterative constraint energy minimizing generalized multiscale finite element method for contact problem}


  \author[CUHK]{Zishang Li\corref{cor1}}
  \cortext[cor1]{Corresponding author.}
  \ead{zsli@math.cuhk.edu.hk}

  \author[CUHK]{Changqing Ye}
  \author[CUHK]{Eric T.~Chung}




  \affiliation[CUHK]{
    organization={Department of Mathematics, The Chinese University of Hong Kong},
    city={Shatin},
    country={Hong Kong SAR}
  }

  \begin{abstract}

    This work presents an Iterative Constraint Energy Minimizing Generalized Multiscale Finite Element Method (ICEM-GMsFEM) for solving the contact problem with high contrast coefficients. The model problem can be characterized by a variational inequality, where we add a penalty term to convert this problem into a non-smooth and non-linear unconstrained minimizing problem. The characterization of the minimizer satisfies the variational form of a mixed Dirilect-Neumann-Robin boundary value problem. So we apply CEM-GMsFEM iteratively and introduce special boundary correctors along with multiscale spaces to achieve an optimal convergence rate. Numerical results are conducted for different highly heterogeneous permeability fields, validating the fast convergence of the CEM-GMsFEM iteration in handling the contact boundary and illustrating the stability of the proposed method with different sets of parameters. We also prove the fast convergence of the proposed iterative CEM-GMsFEM method and provide an error estimate of the multiscale solution under a mild assumption.

  \end{abstract}




\end{frontmatter}


\section{Introduction}
Composite materials have gained significant prominence in various industries and natural settings. The combination of different materials in composites allows for the creation of new materials with enhanced performance characteristics, making them highly sought after in fields such as aerospace, automotive, and construction, among others. 
Problems for composites arising from physics and engineering often exhibit multiple scales and high-contrast features. Traditional methods require very fine grids to solve those problems accurately, which leads to a great number of degrees of freedom and expensive computation. There have been many existing approaches in the literature to handle multiscale problems. These multiscale approaches include multiscale finite element methods \cite{Chen2003,Efendiev2009,Hou1997,Hou1999}, heterogeneous multiscale methods \cite{Abdulle2012,Ming2005,Weinan2007}, variational multiscale methods \cite{Hughes1995,Hughes2007}, generalized finite element methods \cite{Babuska2020,Babuska2011}, generalized multiscale finite element methods (GMsFEM) \cite{Efendiev2013,Chung2014,Chung2016}, and localized orthogonal decomposition methods \cite{Hellman2017,Henning2014,Maalqvist2014,Maalqvist2020}, etc. Most of these approaches are based on encoding fine-scale information into basis functions of finite element methods (FEMs), then solving original problems on multiscale finite element spaces whose dimensions have been greatly reduced compared to default FEMs. Many existing approaches in the field are focused on homogeneous or linear boundary value problems (BVPs) as model problems to study convergence theories and conduct numerical experiments. The extensions of those methods to non-smooth BVPs are relatively not easy to implement due to the presence of discontinuities, singularities, and other forms of non-smooth behavior in the solution or its derivatives. Given the high demand for solving non-smooth BVPs in practical applications, we will examine the effectiveness of current multiscale computational methods and further develop them to handle non-smooth mixed BVPs.

In many industrial applications or engineering problems, contact between deformable bodies plays a crucial role. The Signorini problem is a specific type of unilateral contact boundary value problem in partial differential equations that arises in the study of contact mechanics in solid mechanics. It was first published in the article \cite{Signorini1959} by Antonio Signorini. The solution to the Signorini problem led to the birth of the field of variational inequalities. They were introduced by Fichera in his analysis of the Signorini problem on the elastic equilibrium of a body under unilateral constraints \cite{Fichera1964}. In \cite{Lions1967}, Lions and Stampacchia extended Fichera's analysis to abstract variational inequalities associated with bilinear forms which are coercive or simply nonnegative in real Hilbert spaces as a tool for the study of partial differential elliptic and parabolic equations. Model reduction technique for solving standard PDE has been extended to variational inequality recently, where constructs the reduced basis by combining the greedy algorithm \cite{Rozza2008,Haasdonk2012} or the proper orthogonal decomposition methodology \cite{Fauque2018,Krenciszek2014}. Variational inequality can be characterized by a constrained minimization problem which also provides a numerical treatment \cite{Kikuchi1988}. Despite the fundamental role of contact in the mechanics of solids and structures, contact effects are rarely taken into account in structural analysis. Contact problems are inherently nonlinear \cite{Kikuchi1988}. So the modeling of contact phenomena poses serious difficulties in mathematics and computation. It is far more complex than that encountered in classical linear structural mechanics. In the field of computational contact mechanics, many studies have explored combining contact conditions with discrete formulations and algorithms \cite{Kikuchi1988,Haslinger1981,Kinderlehrer2000,Haslinger1996}. For instance, in \cite{BenBelgacem2003}, the authors revisited three different types of hybrid finite element methods for the Signorini problem. In \cite{Su2022}, the authors used a hybrid variational formulation and applied the contact conditions to the displacement and the stress on the contact zone separately. Despite several previous studies and the rapid improvement in modern computer technology, most finite element software is not fully capable of solving contact problems with robust algorithms \cite{Wriggers2006}. Hence there is still a challenge to design efficient and robust methods for computational contact problems.

\begin{figure}
	\def \svgwidth{\columnwidth}
	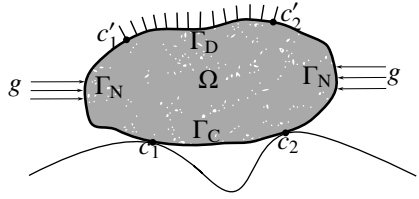
	\caption{A deformable body in contact with a rigid obstacle.}\label{fig:bodyincontact}
\end{figure}

The study and application of the Signorini problem have expanded beyond the physical or mechanical fields. For instance, in hydrostatics, consider a fluid contained in a porous domain $\Omega$ limited partly by a thin membrane $\TC$ which is semi-permeable, meaning that it allows the fluid to pass through only in one direction to get in $\Omega$ \cite{BenBelgacem2003}. In our work, we will focus on heterogeneous medium in mixed contact boundary problems. The representation of this model problem is provided in \cref{fig:bodyincontact}.

The main contribution of our work is as follows. We develop a new iterative multiscale method based on the CEM-GMsFEM for solving the Signorini problem. The variational inequality obtained from the Signorini problem can be characterized by a constrained minimizing problem, where the penalty method is a general technique by adding a penalty term to the objective function that penalizes violations of the constraints. Then it reduces to solving the unconstrained optimization problem by numerical methods. The construction of multiscale space starts with solving a local spectral problem in each coarse element, then the multiscale basis functions are built by solving local constraint energy minimization problems in oversampling domains. The main difficulties in solving the problem arise from the nonlinearity of the contact conditions. Therefore, the semismooth Newton method is introduced during the construction of multiscale space. Specifically, we do not need to update all the basis during iterations, only those near the contact boundary. An analysis of the proposed method is presented. In particular, we show the convergence of our semismooth Newton iteration and give an error estimate. We present numerical results for two different heterogeneous permeability fields and verify that we can get accurate approximations with fewer degrees of freedom using the proposed method. The first case we consider contains small inclusions while the second one contains several high contrast channels. Both examples show our method works well to deal with the contact boundary and the iteration has fast convergence.

We begin, in the next section, by presenting the model of the Signorini problem and the preliminary treatment. Our iterative method and the framework of the multiscale method based on the CEM-GMsFEM are illustrated in Section 3. Then in section 4, we state the error analysis and show the convergence for our method. The numerical experiments are given in Section 5, where we also verify the effectiveness and efficiency of our method. Finally, some conclusions are provided in Section 6.

\section{Preliminaries}\label{sec:pre}
\subsection{Model problem}
We will consider the contact problem of Signorini type for a second-order elliptic partial differential equation. We denote a Lipschitz domain by $\Omega \subset \R^d$ ($d=2\ \text{or}\ 3$), and $\kappa\in L^\infty\left(\Omega ;\ \R^{d \times d} \right) $ a matrix-valued function defined on $ \Omega $ represents a heterogeneous permeability field with high contrast:
\begin{equation}\label{eq:strong}
	\left\{
	\begin{aligned}
		 & - \Div\left( {\kappa \nabla u} \right) = f                                                     &  & \text { in } \Omega, \\
		 & u=0                                                                                            &  & \text { on } \TD,    \\
		 & \kappa\nabla u\cdot\bm{n}=p                                                                    &  & \text { on } \TN,    \\
		 & u \leqslant0, \quad \kappa\nabla u\cdot\bm{n} \leqslant0,\quad (\kappa\nabla u\cdot\bm{n})u =0 &  & \text { on } \TC,    \\
	\end{aligned}
	\right.
\end{equation}
where $\bm{n}$ is the outward unit normal to $\partial\Omega$, $\TD$, $\TN$ and $\TC$ are three nonempty disjointed parts of $\partial\Omega$. In this paper, we present the following assumptions:
\begin{enumerate}
	\item There exist two positive constants $ \kappa' $ and $ \kappa'' $ such that  $0<\kappa' \leqslant \kappa(\bm{x}) \leqslant \kappa''<\infty$ for almost all $ \bm{x}\in\Omega $.
	\item The source term $f \in L^2(\Omega ) $, and the inhomogenous Neumann boundary term $ p \in L^2(\TN) $.
\end{enumerate}

To solve this problem, it is appropriate to use a functional framework that involves a subset of the Sobolev space $H^1(\Omega)$ defined as
\[ V\coloneqq\left\{v \in H^1(\Omega): \  v=0\text{ on } \TD \right\}. \]
The contact condition is then explicitly incorporated in the following closed convex set
\[ K\coloneqq\left\{v \in V: \  v\leqslant0 \text{ on } \TC\right\}.\]
By the primal variational principle for the Signorini problem, the exact solution of the above contact problem (\cref{eq:strong}) is characterized by the variational inequality: find $u \in K$ such that
\begin{equation}\label{eq:var inequ}
	a(u, v-u) \geqslant L(v-u), \  \forall v \in K,
\end{equation}
where
\[  a(u, v) =\int_{\Omega}\kappa \nabla u \cdot\nabla v \di \bm{x} \ \text{and}\ L(v) =\int_{\Omega} f v \di \bm{x}+\int_{\Gamma_N} p v \di \sigma. \]
We introduce the notation for the energy norm $ {\left\| v \right\|_a} \coloneqq \sqrt {a(v,v)}$ on $ \Omega $. For a subdomain $ \omega \subset \Omega $ , we also introduce the norm $ {\left\| v \right\|_{a(\omega)}} \coloneqq \sqrt{\int_{\omega}\kappa \nabla u \cdot\nabla v \di \bm{x}}. $

By Theorem 3.9 (see \cite{Kikuchi1988}), the weak problem (\cref{eq:var inequ}) is well-posed and the unique solution of this variational inequality exists. Moreover, this solution is also the minimizer of the constrained minimization problem: given the functional $F: V \rightarrow \R$ is of the form
\[F(v)=\frac{1}{2} a(v, v)-L(v),\]
find $u \in K$ such that
\begin{equation}\label{eq:con min}
	F(u) = \inf\limits_{v \in K} F(v).
\end{equation}

\subsection{Penalty method}\label{subsec:penalty}
Constrained minimization problems involve optimizing functions subject to constraints and are generally more complex than unconstrained problems. Meanwhile, the unconstrained minimization problems have more well-established solving methods and find extensive applications in various industries and engineering fields, where they minimize a function without any constraints on the variables. So the primary objective of our research is to transform constrained minimization problems into unconstrained problems, and one commonly employed approach for achieving this conversion is the penalty method. The penalty method involves the introduction of additional penalty terms to the objective function, which effectively penalizes violations of the constraints. By incorporating these penalty terms, the original constrained problem is transformed into an unconstrained problem, allowing for the application of unconstrained optimization techniques and algorithms. Here, we introduce the penalty term $P:V\rightarrow \R$ and a new functional $F_\varepsilon:V\rightarrow\R $ depending on a real parameter $\varepsilon >0$, of the form
\[ {F_\varepsilon }(v) = F(v) + \frac{1}{\varepsilon }P(v). \]
Hence, in accordance with Theorem 3.1 in \cite{Kikuchi1988}, for each $\varepsilon>0$ there exists a ${u_\varepsilon } \in V$ which minimizes $F_\varepsilon$. The corresponding constrained minimization problem is: find $u_\varepsilon \in V$ such that
\begin{equation}\label{eq:con min penalty}
	{F_\varepsilon }({u_\varepsilon })= \inf\limits_{v \in V} {F_\varepsilon }(v).
\end{equation}
Moreover, if the functionals $F$ and $P$ are Gateaux-differentiable, the minimizer ${u_\varepsilon }$ can be characterized by
\[ \left\langle {DF(u_\varepsilon),v} \right\rangle +\frac{1}{\varepsilon}\left\langle {DP({u_\varepsilon }),v} \right\rangle = 0, \]
where $ \left\langle {\cdot,\cdot} \right\rangle $ is the functional dual on $ V $. To apply the penalty method to the Signorini condition, the penalty functional $P: V \rightarrow \R$ is taken as
\[ P(v) = {\frac{1}{2}}\int_{\TC } {v_ + ^2\di\sigma }, \]
which satisfies a property: the more a candidate minimizer $v \in V$ violates the constraint $v\leqslant0$ on the contact boundary $\TC$, the greater the penalty that must be paid. This functional $P$ is also Gateaux-differentiable on $V$ and
\[ \begin{aligned}
		\left\langle {DP(u),v} \right\rangle
		 & =\mathop {\lim }\limits_{t \to 0} \frac{P(u+tv)-P(u)}{t}                                                               \\
		 & =\mathop {\lim }\limits_{t \to 0^+} \frac{1}{2t}\int_{\TC} {{(u+tv)}_ + ^2-{u}_ + ^2\di\sigma }                        \\
		 & =\mathop {\lim }\limits_{t \to 0^+} \frac{1}{2t}\int_{\TC } {\int_0^1 {\frac{\di}{\di s}(u+stv)_+^2 \di s} \di\sigma } \\
		 & =\mathop {\lim }\limits_{t \to 0^+} \frac{1}{2t}\int_{\TC } {\int_0^1{2(u+stv)_+tv \di s} \di\sigma }                  \\
		 & =\int_{\TC } {{u_ + }v\di\sigma }.                                                                                     \\
	\end{aligned} \]
Thus, the characterization of the minimizer of the unconstrained minimization problem is: find $u_\varepsilon \in V$, such that
\begin{equation}\label{eq:char of the minimizer of F_epsilon}
	a(u_\varepsilon, v)+\frac{1}{\varepsilon}\int_{\TC } {{(u_\varepsilon)_ + }v\di\sigma }=L(v),\ \forall v\in V.
\end{equation}

\subsection{Semismooth Newton Method}\label{subsec:semi-new}
To find the minimizer of the unconstrained minimization problem, we take the residual functional for the characterization \cref{eq:char of the minimizer of F_epsilon} as
\begin{equation}\label{eq:R(u)}
	R(u) =a(u, v)+\frac{1}{\varepsilon}\int_{\TC } {g(u)v\di\sigma }-L(v),
\end{equation}
where we denote $ g(u) = {u_ + } = \max \left\{ {0,\ u} \right\} $.
Then we need to find the root $u_\varepsilon \in V$, such that $R(u_\varepsilon)=0$.
Since $ R(u) $ lacks the necessary differentiability in the classical sense, we must turn to generalized Newton methods for this nonlinear problem. One such generalized method is the semismooth Newton method. The definition of semismooth is cumbersome to handle \cite{Hintermueller2010}. To simplify, we introduce the following theorem which provides equivalent characterizations to our case.
\begin{theorem}
	Let $ R: V\to \R $. Then, for $ u\in V$, the following statements are equivalent:\\
	(a) $ R $ is semismooth at $ u $.\\
	(b) $ R $ is locally Lipschitz continuous at $ u $, $  R'(u;\cdot) $ exists, and for any $ G\in \partial R(u + d) $,
	\[\norm{Gd-R'(u, d)} = \bigO(\norm{d})\text{ as } d\to 0 .\]\\
	(c) $ R $ is locally Lipschitz continuous at $ u $, $  R'(u;\cdot) $ exists, and for any $ G\in \partial R(u + d) $,
	\[\norm{R(u+d)-R(u)-Gd}=\bigO(\norm{d}) \text{ as } d\to 0 .\]
\end{theorem}
Here $ R'(u; d) $ and $ \partial R(u + d) $ are the directional derivative and collection of the generalized directional derivative of $ R $ at $ u $ in direction $ d $ respectively. The definition can be found in \cite{Hintermueller2010}. We find that $ R(u) $ is semismooth on $ V $. Then a generalized derivative is defined for the semismooth Newton procedure.

\begin{definition}
	The mapping $ R:V \to \R $ is Newton differentiable on the open set $ U \subset V$, if there exists a family of mappings $ G  : U \to L(V,V')$ such that
	\[\lim\limits_{d \to 0}\frac{1}{\norm{d}} \norm{R(u+d)-R(u)-G(u+d)d} =0. \]
	for every $ u \in U $. The operator $ G $ is referred to as Newton's derivative of $ R $.
\end{definition}

Note that $ R(u) $ is Newton differentiable on $ V $. The Newton derivative of $ R(u) $ is
\[\left\langle {G(u)d,v} \right\rangle \coloneqq a(d,v) + \frac{1}{\varepsilon }\int_{\TC }g'({u})dv\di\sigma,\  \forall d,v\in V,\]
where the derivative of $ g(u) $ is
\[ g^{\prime}(u)=\mathbbm{1}(u>0)=\left\{\begin{aligned}
		 & 1,\ u>0,           \\
		 & 0,\ u \leqslant 0.
	\end{aligned}\right. \]

Similar to the classical Newton method, the semismooth Newton algorithm is stated in \cref{algo}:\\
\begin{algorithm}
	\caption{Semismooth Newton method} \label{algo}
	\hspace*{0.02in} {\bf Input:} 
	$ R $, $ u_0 $, $\mathup{tol} $
	\begin{algorithmic}[1]
		\State $ R$ is Newton differentiable in $V$, $ u_0 \in V$, $ k=0 $
		\While{$ k<1 $ or $ R(u_{k+1})>\mathup{tol} $}
		\State Solve $ G(u_k)d_k=R(u_k) $ to obtain $ d_k $
		\State $ u_{k+1} = u_k-d_k $, $ k = k+1 $
		\EndWhile
		\State \Return $u_{k+1}$
	\end{algorithmic}
	\hspace*{0.02in} {\bf Output:} 
	the numerical solution $u_{k+1}$
\end{algorithm}\\

The main procedure in variational form is: find $ u_k $, such that
\begin{equation}\label{eq:main procedure in Newton}
	\left\{
	\begin{aligned}
		 & a({d_k},v) + \frac{1}{\varepsilon }\int_{\TC } g'({u_k}){d_k}v\di\sigma= a(u_k, v)+\frac{1}{\varepsilon}\int_{\TC } {g({u_k })v\di\sigma }-L(v),\  \forall v\in V, \\
		 & {u_{k + 1}} = {u_k} - {d_k}.
	\end{aligned}
	\right.
\end{equation}

\section{Numerical Method}\label{sec:method}
\subsection{Conversion of the contact problem}\label{subsec:conversion}
Although many theoretical derivations and problem transformations have been done for the Signorini contact problem. We still need to do some processing before using numerical methods and getting the numerical solution. The characterization of the minimizer of $F_\varepsilon$ \cref{eq:char of the minimizer of F_epsilon} in \cref{subsec:penalty} can be considered as the weak form of the specific partial differential equation. Remarkably, it is identical with the weak form of the following Neumann boundary problem
\begin{equation}\label{eq:Neumann}
	\left\{
	\begin{aligned}
		 & - \Div\left( {\kappa \nabla u} \right) = f,                         &  & \text { in } \Omega, \\
		 & u=0,                                                                &  & \text { on } \TD,    \\
		 & \kappa\nabla u\cdot\bm{n}=p,                                        &  & \text { on } \TN,    \\
		 & \kappa\nabla u\cdot\bm{n}=-\frac{1}{\varepsilon }g\left( u \right), &  & \text { on } \TC.    \\
	\end{aligned}
	\right.
\end{equation}
The weak formulation of \cref{eq:Neumann} is to find $u \in V$ such that
\[a(u, v)=f(v)-\frac{1}{\varepsilon}\int_{\TC } {{u_ + }v\di\sigma }+\int_{\TN}p v\di\sigma, \quad\forall v\in V,\]
where
\[ a\left( {w,v} \right) = \int_\Omega {\kappa \nabla w\cdot \nabla v} \di\bm{x}\  \text{and}\ f(v)=\int_{\Omega} f v \di \bm{x}. \]
Consequently, it is possible that the numerical methods for solving Neumann boundary value problems can be applied to determine the contact problem minimizer.

To linearize the function $g(u)$, we can utilize the first-order Taylor series approximation. This approximation expresses $g(u)$ as follows:
\[g\left( u \right) \approx g\left( {{u_k}} \right) + g'\left( {{u_k}} \right)\left( {u - {u_k}} \right),\]
where $g'(u_k)$ represents the derivative of $g(u)$ evaluated at the point $u_k$. Thus, we have
\[\frac{1}{\varepsilon }g'\left( {{u_k}} \right)u + \kappa\nabla u\cdot\bm{n} = - \frac{1}{\varepsilon }\left[ {g\left( {{u_k}} \right) - g'\left( {{u_k}} \right){u_k}} \right].\]
Take ${b_k}$ and ${q_k}$ as
\[ {b_k} = \frac{1}{\varepsilon }g'\left( {{u_k}} \right)\text{ and }{q_k} = - \frac{1}{\varepsilon }\left[ {g\left( {{u_k}} \right) - g'\left( {{u_k}} \right){u_k}} \right].\]
Then \cref{eq:Neumann} can be converted into an inhomogeneous Robin boundary value problem
\begin{equation}\label{eq:newcon}
	\left\{
	\begin{aligned}
		 & -\Div\left( {\kappa \nabla u} \right) = f, &  & \text { in } \Omega, \\
		 & u=0,                                       &  & \text { on } \TD,    \\
		 & \kappa\nabla u\cdot\bm{n}=p,               &  & \text { on } \TN,    \\
		 & b_ku+\kappa\nabla u\cdot\bm{n}=q_k,        &  & \text { on } \TC.    \\
	\end{aligned}
	\right.
\end{equation}
We rewrite this mixed boundary problem \cref{eq:newcon} in a variational form
\[ \int_\Omega {\kappa \nabla u\cdot\nabla v} \di\bm{x} +\int_{\TC } {b_kuv} \di\sigma=\int_\Omega fv \di\bm{x}+\int_{\TN}p v\di\sigma +\int_{\TC} q_k v\di\sigma, \]
and the bilinear form satisfies
\[ \tilde{a}\left( {w,v;b_k} \right) = \int_\Omega {\kappa \nabla w\cdot\nabla v} \di\bm{x} + \int_{\TC } {{b_k}wv} \di\sigma =a\left( {w,v} \right)+ \int_{\TC } {{b_k}wv} \di\sigma. \]
Here we introduce a notation for the bilinear form $ \tilde{a}\left( {\cdot,\ \cdot\ ;\ b} \right)$ :
\[ \tilde{a}\left( {w,v;b} \right) \coloneqq \int_\Omega {\kappa \nabla w\cdot\nabla v} \di\bm{x} + \int_{\TC } {{b}wv} \di\sigma. \]
Then we have
\begin{equation}\label{var_robin}
	a\left( {u,v} \right)+ \frac{1}{\varepsilon }\int_{\TC } {g'\left( {{u_k}} \right)uv} \di\sigma=\int_\Omega fv \di\bm{x}+\int_{\TN}p v\di\sigma - \frac{1}{\varepsilon }\int_{\TC} {{ \left[ {g\left( {{u_k}} \right) - g'\left( {{u_k}} \right){u_k}} \right]}} v\di\sigma.
\end{equation}
This variational equation is exactly equivalent to the \cref{eq:main procedure in Newton} obtained by the semismooth Newton iteration in \cref{subsec:semi-new} if we take $u=u_k-d_k$. In the case of Signorini typed contact problem, the function $g(u) = (u)_+$. Then we have
\[ b_k\left(u_k\right)=\frac{1}{\varepsilon} g^{\prime}\left(u_k\right)=\frac{1}{\varepsilon} \mathbbm{1}\left(u_k>0\right)\]
and

\[ \begin{aligned}
		{q_k}\left( {{u_k}} \right) & = - \frac{1}{\varepsilon }\left[ {g\left( {{u_k}} \right) - g'\left( {{u_k}} \right){u_k}} \right] \\
		                            & =- \frac{1}{\varepsilon }\left[ (u_k)_+ -\mathbbm{1}(u_k>0) (u_k)_+ \right]                        \\
		                            & =- \frac{1}{\varepsilon }\left[ (u_k)_+ - (u_k)_+ \right]                                          \\
		                            & =0.
	\end{aligned} \]
Thus, the variational form \cref{var_robin} can be simplified as: find $u_{k+1}\in V$, such that
\begin{equation}\label{eq:var_iter}
	a(u_{k+1}, v)+\frac{1}{\varepsilon}\int_{\TC} {u_{k+1}v\cdot\mathbbm{1}(u_k>0)\di\sigma }=\int_\Omega fv \di\bm{x}+\int_{\TN}p v\di\sigma,\  \forall v\in V.
\end{equation}
By doing the iteration, we can obtain the minimizer of the unconstrained functional $F_\varepsilon$. Drawing upon the previous analysis, it is feasible to transform the contact problem into a mixed Robin boundary problem. This transformation holds significant importance as it may allow us to effectively handle this problem by modifying the established finite element tools and methodologies.

\subsection{One step CEM-GMsFEM solver}\label{subsec:cemgmsfem}
\begin{figure}
	\centering
	\begin{tikzpicture}[scale=1.2]
		\draw[step=0.25, gray, thin] (-0.4, -0.4) grid (5.4, 5.4);
		\draw[step=1.0, black, very thick] (-0.4, -0.4) grid (5.4, 5.4);
		\foreach \x in {0,...,5}
		\foreach \y in {0,...,5}{
				\fill (1.0 * \x, 1.0 * \y) circle (1.5pt);
			}

		\fill[MyColor4, opacity=0.6] (0.0, 0.0) rectangle (5.0, 5.0);
		\fill[MyColor1, opacity=0.6] (2.0, 2.0) rectangle (3.0, 3.0);
		\fill[MyColor2, opacity=0.6] (3.5, 0.25) rectangle (3.75, 0.5);

		\node at (2.5, 2.5) {$K_i$};
		\node at (3.5, 3.5) {$K^2_i$};
		\node at (3.625, 0.375) {$\tau$};

		\fill[MyColor2, opacity=0.6] (5.5, 2.4) rectangle (5.9, 2.6);
		\node[right] at (5.9, 2.5) {=\ a fine element $\tau$};

		\fill[MyColor1, opacity=0.6] (5.5, 1.9) rectangle (5.9, 2.1);
		\node[right] at (5.9, 2.0) {=\ a coarse element $K_i$};

		\fill[MyColor4, opacity=0.6] (5.5, 1.4) rectangle (5.9, 1.6);
		\node[right] at (5.9, 1.5) {=\ an oversampled domain $K^m_i$ (here $ m=2 $)};
	\end{tikzpicture}
	\caption{Illustration of meshes, fine element, coarse element, and oversampling domain.}
	\label{fig:grid}
\end{figure}
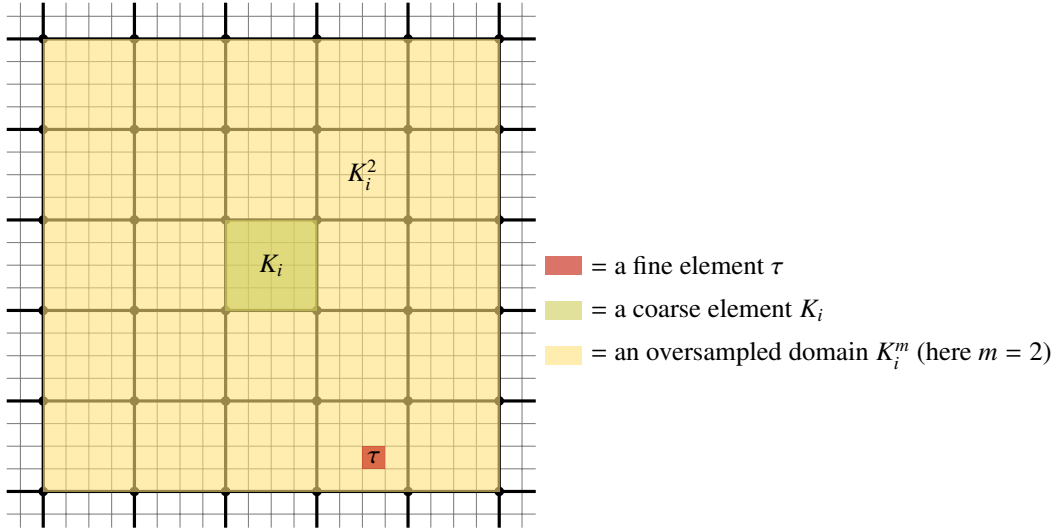

We will focus on the numerical solution of the variational form \cref{eq:var_iter} in one iteration by the CEM-GMsFEM method. Consider a conforming partition $\mathcal{T}^H$ of a domain $\Omega$ into $ N $ finite elements $K_i$, where $H$ denotes the coarse-mesh size. This is to be distinguished from another fine mesh $\mathcal{T}^h$ with the mesh size $ h $ and will be utilized to compute multiscale basis functions. The fine element is represented as $ \tau $. For each coarse element $K_i \in \mathcal{T}^H$ with $1 \leqslant i \leqslant N$, we define an oversampling domain $K^m_i(m \geqslant 1)$ as the domain obtained by augmenting the coarse element $K_i$ with an additional $m$ layers of neighboring coarse elements. A representation of the fine grid, coarse grid, and oversampling domain is provided in \cref{fig:grid}.
For this quadrilateral mesh, there are 4 vertices contained in an element. We can construct a set of Lagrange bases $  \left\{ {\eta _i^1,\eta _i^2, \eta _i^3 ,\eta _i^4} \right\} $ of the coarse element $K_i \in \mathcal{T}^H$. Then we define $ \tilde \kappa (x) $ piecewise by
\begin{equation}\label{eq:kappa}
	\tilde \kappa (x)=3\sum\limits_{j = 1}^4 {\kappa (x)\nabla \eta _i^j}\cdot \nabla \eta _i^j
\end{equation}
in $ K_i $ which will be used in the following spectral problem.

The process of the construction of CEM-GMsFEM basis functions can be divided into two stages. The first stage involves the construction of the auxiliary space by solving a local spectral problem in each coarse element $ K_i$: find eigen-pairs $ \{ {\lambda _i^j,\phi _i^j} \} $ such that
\begin{equation}\label{eq:eigen_pro}
	{\tilde{a}_i}(\phi _i^j,v;b) = \lambda _i^j{s_i}(\phi _i^j,v), \  \forall v\in H^1(K_i).
\end{equation}
where
\[ \tilde{a}_i(u, v; b) =\int_{K_i}\kappa \nabla u \cdot\nabla v \di \bm{x}+ \int_{\TC \cap \partial {K_i}} {buv}\di\sigma \ \text{and}\ s_i(u,v)= \int_{{K_i}} \tilde \kappa uv\di\bm{x}. \]
The bilinear form $s(w,v) \coloneqq \int_{{\Omega}} \tilde \kappa wv\di\bm{x} $, and note that $ s(w, v) $ can be well-defined on $ L^2(\Omega ) $. Similarly, we denote the norm
\[ {\left\| v \right\|_s} \coloneqq \sqrt {s(v,v)}\ \text{and}\ {\left\| v \right\|_{s(K_i)}} \coloneqq \sqrt{s_i(v,v)}. \]
Let the eigenvalues $ \{\lambda_i^j\}_{j=1}^\infty $ be arranged in an ascending order, we notice that $\lambda_i^1 = 0 $ always holds. We define the local auxiliary multiscale space $ V_i^{\mathup{aux}} $ by using the first $ l_i $ eigenfunctions
\[V_i^{\mathup{aux}} \coloneqq \Span \CurlyBrackets*{\phi_i^j:1 \leqslant j \leqslant l_i },\]
and the orthogonal projection $ \pi_i :V(K_i)\rightarrow V_i^{\mathup{aux}} $ with respect to the inner product $ s(\cdot , \cdot ) $ is
\[\pi_i(v)\coloneqq\sum\limits_{j = 1}^{l_i}\frac{s(\phi_i^j,v)}{s(\phi_i^j,\phi_i^j)}\phi_i^j.\]
Then the global auxiliary space $V^{\mathup{aux}} $ is defined by using these local auxiliary spaces
$ V^{\mathup{aux}}=\mathop  \oplus_{i = i}^N V_i^{\mathup{aux}}, $
and the global projection is $ \pi \coloneqq \sum^N_{i=1} \pi_i $ accordingly.

In the second stage, the multiscale basis functions are formed by solving some local constraint energy minimization problems in oversampling domains $ K_i^m$. Let
\[ V(K_i^m)=\left\{v \in H^1(K_i^m): \  v=0\text{ on } \TD\cap {\partial K_i^m}\ \text{or}\ \Omega\cap{\partial K_i^m} \right\}, \]
for each oversampling domain $ K_i^m $. Then, the multiscale basis functions are defined by
\begin{equation}\label{eq:ms_min}
	\psi_i^{j,m} = \mathup{argmin}\left\{ \tilde{a}(\psi , \psi; b ) + s(\pi \psi  - \phi^j_i , \pi \psi  - \phi^j_i ):\ \psi \in V(K_i^m) \right\},
\end{equation}
which is a relaxed version of the energy minimization problem \cite{Chung2018}
\[\psi_i^{j,m}= \mathup{argmin}\left\{ \tilde{a}(\psi,\psi;b): \psi \in V(K_i^m),\  \pi \psi=\phi^j_i \right\}.\]
We note that \cref{eq:ms_min} is equivalent to: find $ \psi_i^{j,m} $ such that
\begin{equation}\label{eq:char_ms_min}
	\tilde{a}(\psi^{j,m}_i , v;b) + s(\pi \psi^{j,m}_i , \pi v) = s(\phi^j_i , \pi v),\  \forall v \in  V(K_i^m).
\end{equation}
Then the CEM-GMsFEM multiscale finite element space is defined by
\[ V_\mathup{ms}^m=\mathup{span}\left\{\psi_i^{j,m}:1 \leqslant j \leqslant l_i, 1 \leqslant i \leqslant N \right\}. \]

Given that the error analysis theory of the original CEM-GMsFEM strongly relies on the existence of $ L^2 $ source term, we might consider introducing function ${\cal N}_i^m{p}$ by imitating the construction of multiscale bases. This approach aims to extend the applicability of the CEM-GMsFEM to cases where an $L^2$ source term is not readily available or suitable. The function ${\cal N}_i^m{p}$ is defined by solving the following local problem (\cref{eq:corrector}) in the oversampling domain $K_i^m$. We can see that the zero-extension of a function in $ V(K_i^m) $ still belongs to $ V $.

Denote this process in one iteration as the  CEM-GMsFEM solver. The output numerical solution of this solver is $ {\cal S}\{b(\hat u),\kappa,f,p\} $, where $ \hat u $ is the substitutable input function. When the input function ${\hat u}$ is given, it is easy to compute ${b(\hat u)}$. We will perform the following four steps to obtain the new multiscale numerical solution:
\begin{description}
	\item[STEP 1] Find ${\cal N}_i^m{p} \in V(K_i^m)$ such that,
	      \begin{equation}\label{eq:corrector}
		      \tilde{a}\left( {{\cal N}_i^m{p},v;b(\hat u)} \right) + s\left( {\pi {\cal N}_i^m{p},\pi v} \right) = \int_{\partial {K_i} \cap \TN } {pv}\di\sigma,\  \forall v \in V(K_i^m),
	      \end{equation}
	      then obtain
	      \[ {{\cal N}^m}{p} = \sum\limits_{i = 1}^N {{\cal N}_i^m} {p}. \]
	\item[STEP 2] Construct the auxiliary space $V^{\mathup{aux}} $ by \cref{eq:eigen_pro} and the multiscale function space $V_{\mathup{ms}}^m$ by \cref{eq:char_ms_min}.
	\item[STEP 3] Solve ${w^m} \in V_{\mathup{ms}}^m$ such that for all $v \in V_{\mathup{ms}}^m$,
	      \[\tilde{a}\left( {{w^m},v;b(\hat u)} \right) = \int_\Omega f v\di\bm{x} + \int_{{\Gamma _{\rm{N}}}} {{p}} v\di\sigma - \tilde{a}\left( {{{\cal N}^m}{p},v;b_k} \right).\]
	\item[STEP 4] Construct the numerical solution to approximate the real solution as
	      \[{\cal S}\{b(\hat u),\kappa,f,p\} \approx {w^m} + {{\cal N}^m}{p}.\]
\end{description}

\subsection{Iterative CEM-GMsFEM Algorithm}\label{subsec:itercemgmsfem}
The iterative CEM-GMsFEM method involves using the CEM-GMsFEM solver iteratively to solve the contact problem. For the following iteration, the Robin coefficient $b$ and new numerical solution $ \cal S $  can be computed by inputting the obtained numerical solution which involves executing the procedures in \cref{subsec:cemgmsfem} to derive the updated multiscale solution. The condition for iteration termination is set as the error between the numerical solutions $ u_{k-1}^\mathup{cem} $ and $ u_k^\mathup{cem} $ is less than or equal to the tolerance (TOL):
\begin{equation}\label{eq:tolerance}
	\norm{u_{k}^\mathup{cem}-u_{k-1}^\mathup{cem}}\leqslant\mathup{tol}.
\end{equation}
We will conduct iterations until it satisfies this termination condition. The proposed multiscale method is summarized as in \cref{algo:cem}.
\begin{algorithm}
	\caption{Iterative CEM-GMsFEM for contact problem} \label{algo:cem}
	\hspace*{0.02in} {\bf Input:} 
	coefficient $\kappa$, source term $ f $, Neumann boundary term $ p $, $\mathup{tol} $
	\begin{algorithmic}[1]
		\State $ k=0 $ 
		\While{$k<1$ or $ \norm{u_{k}^\mathup{cem}-u_{k-1}^\mathup{cem}}>\mathup{tol} $} 
		\State $ k=k+1 $
		\State Use CEM-GMsFEM solver $ {\cal S}\{b(u_{k-1}^\mathup{cem}),\kappa,f,p\}$ to obtain $u_{k}^\mathup{cem}$
		\EndWhile
		\State \Return $ k $ and $u_{k}^\mathup{cem}$
	\end{algorithmic}
	\hspace*{0.02in} {\bf Output:} 
	the iterative number $ k $ and the numerical solution $u_{k}^\mathup{cem}$
\end{algorithm}

\section{Numerical experiments}\label{sec:num}
\begin{figure}
	\centering
	\begin{subfigure}[b]{0.312\textwidth}
		\centering
		\includegraphics[width=\textwidth]{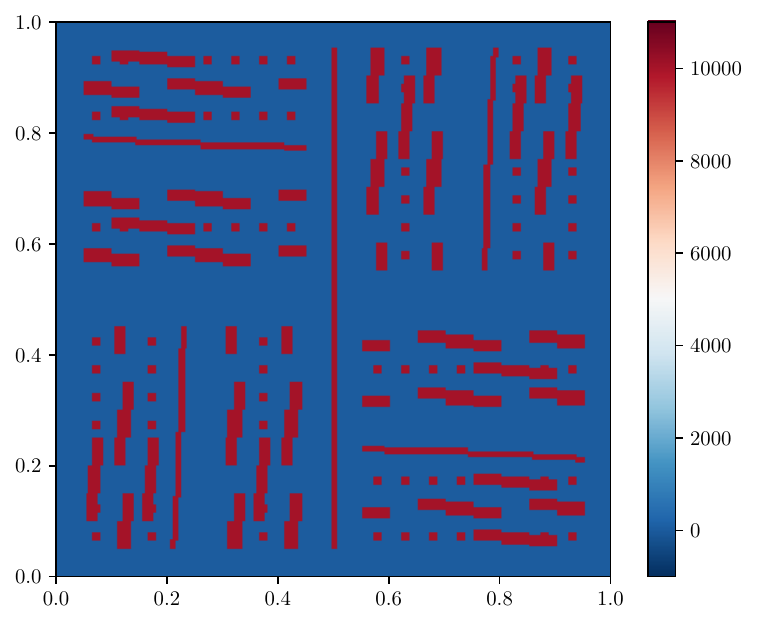}
		\caption{}\label{fig:Medium A}
	\end{subfigure}
	\hfill
	\begin{subfigure}[b]{0.312\textwidth}
		\centering
		\includegraphics[width=\textwidth]{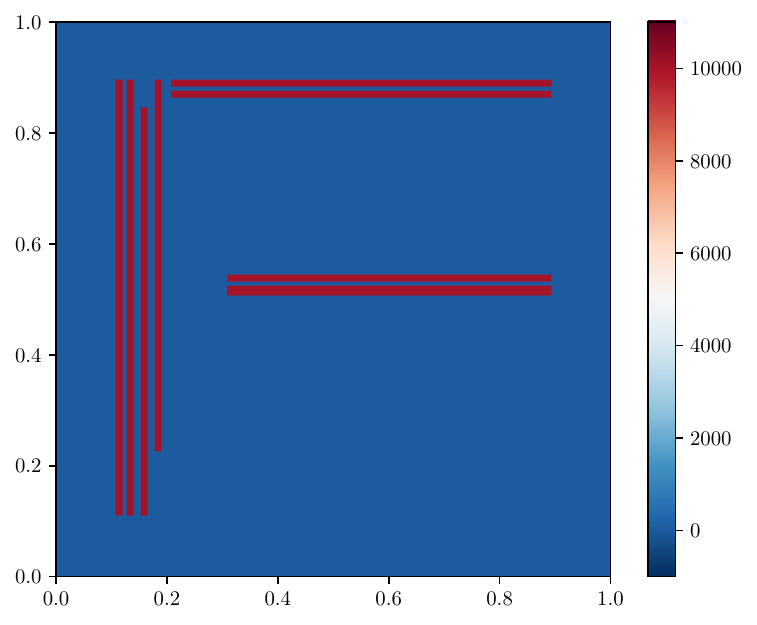}
		\caption{}
		\label{fig:Medium B}
	\end{subfigure}
	\hfill
	\begin{subfigure}[b]{0.312\textwidth}
		\centering
		\includegraphics[width=\textwidth]{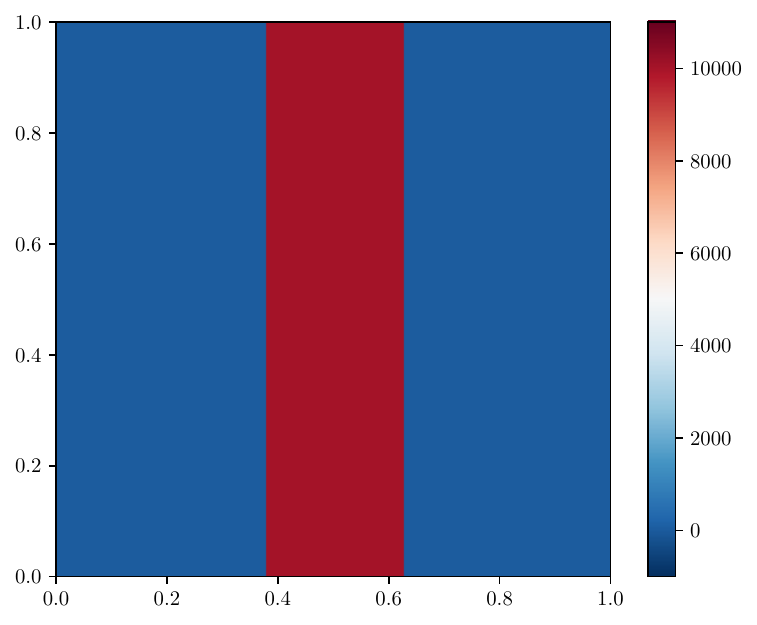}
		\caption{}
		\label{fig:Medium C}
	\end{subfigure}
	\caption{The permeability fields (a) medium A; (b) medium B; (c) medium C.}\label{fig:medium}
\end{figure}

\begin{figure}
	\centering
	\begin{subfigure}[b]{0.3\textwidth}
		\centering
		\includegraphics[width=\textwidth]{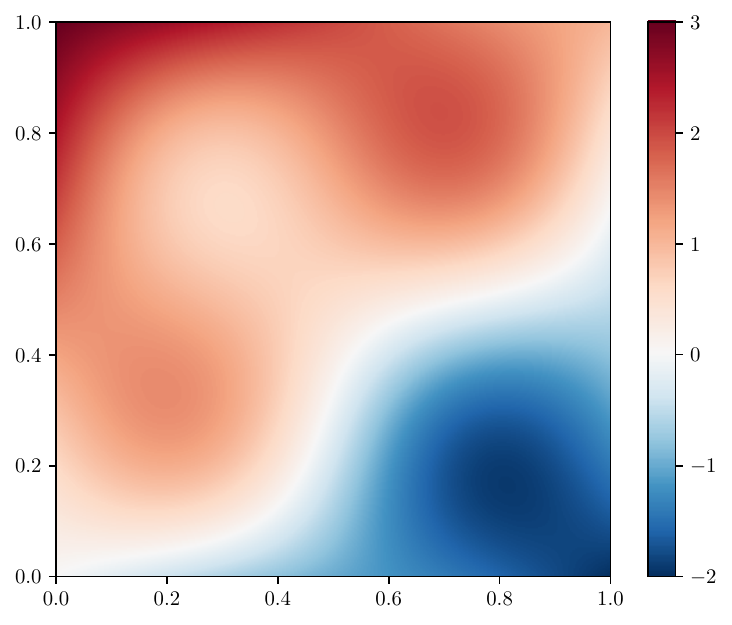}
		\caption{}\label{fig:f_1}
	\end{subfigure}
	\hfill
	\begin{subfigure}[b]{0.312\textwidth}
		\centering
		\includegraphics[width=\textwidth]{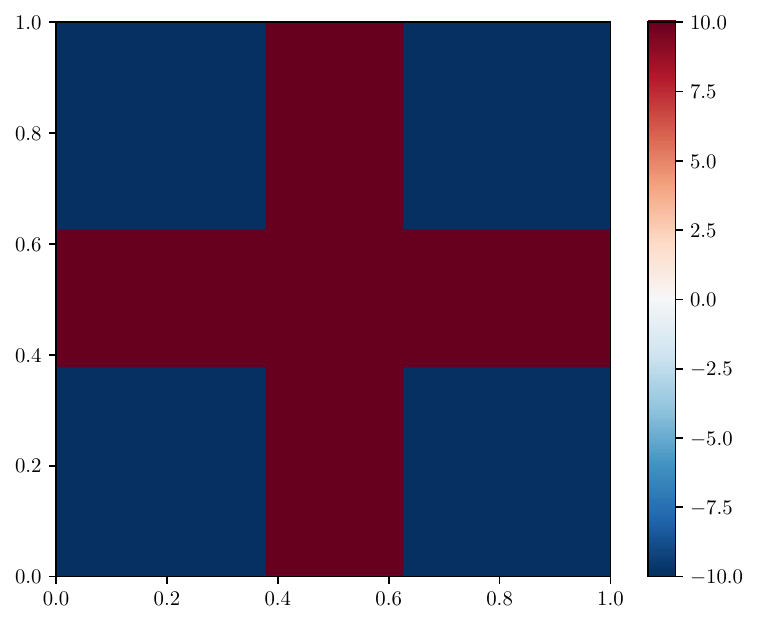}
		\caption{}\label{fig:f_2}
	\end{subfigure}
	\hfill
	\begin{subfigure}[b]{0.312\textwidth}
		\centering
		\includegraphics[width=\textwidth]{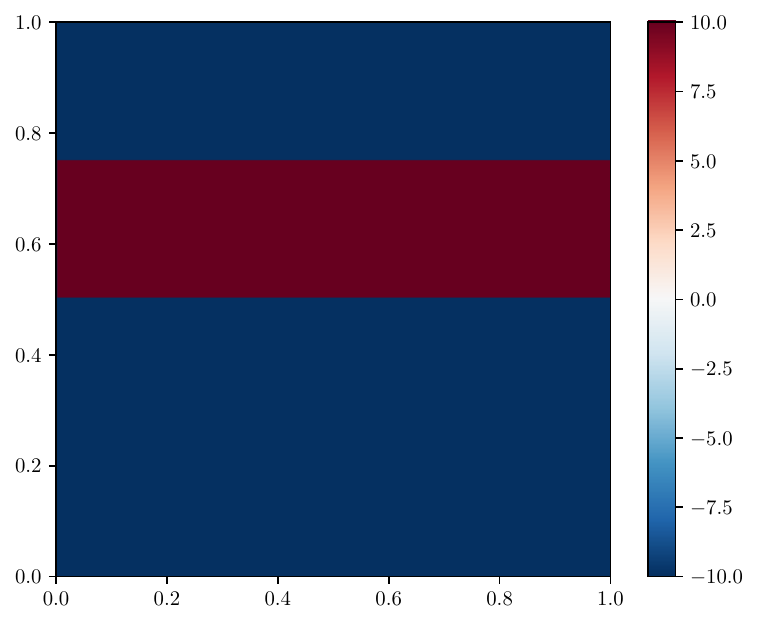}
		\caption{}\label{fig:f_3}
	\end{subfigure}
	\caption{The source function (a) $ f_1 $; (b) $ f_2 $; (c) $ f_3 $.}\label{fig:f}
\end{figure}

In this section, we will present some numerical experiments to demonstrate that the multiscale method proposed is effective in a high contrast coefficient setting. We set the domain $ \Omega = (0, 1) ^2 $. The medium parameter $ \kappa $ has a $ 400 \times 400 $ resolution and only takes two values. The matrix phase value is $\kappa_\mathup{m}=1 $ and the value in the channels and inclusions is $\kappa_\mathup{I} \gg \kappa_\mathup{m}$, as shown in \cref{fig:medium}.
We define contrast ratios $ \kappa_\mathup{R} \coloneqq\kappa_\mathup{I}/\kappa_\mathup{m}$. For simplifying the implementation, we especially choose $ \tilde \kappa = 24\kappa /H^2 $ instead of the original definition in \cref{eq:kappa}. Moreover, we always fix the number of eigenvectors used to construct auxiliary space $ V^{\mathup{aux}}_i $ as $l_\mathup{m}$, i.e., $ l_1 = l_2 = \cdot \cdot \cdot = l_N = l_\mathup{m} $. The source functions $ f $ are given by $f_1(x,y)=-2x+3y+\sin (2\pi x)\sin (2\pi y)$ for any $ (x,y)\in\Omega $,
\begin{equation}\label{f2}
	f_2=\left\{ \begin{aligned}
		10,  & \ 0 < x < 1\text{ and }\frac{3}{8}< y < \frac{5}{8}, \\
		10,  & \ \frac{3}{8}< x < \frac{5}{8}\text{ and }0 < y< 1,  \\
		-10, & \text{ else},
	\end{aligned} \right.
\end{equation}
and
\begin{equation}\label{f3}
	f_3=\left\{ \begin{aligned}
		10,  & \ 0 < x < 1\text{ and }\frac{1}{2}< y < \frac{3}{4}, \\
		-10, & \ \text{ else},
	\end{aligned} \right.
\end{equation}
in \cref{fig:f}.
The reference solutions $ u^{\mathup{fe}}$ are obtained by the bilinear Lagrange FEM with a fine mesh $ 400 \times 400 $. The default parameters are taken as follows: coarse mesh sizes $ H =1/80 $, penalty parameter $ \varepsilon=10^{-4} $, contrast ratios $ \kappa_\mathup{R} = 10^3$, eigenvector numbers $ l_\mathup{m}= 4 $ and oversampling layers $ m=4 $.

After the $ k $-th iteration, the numerical solution $ u_{k}^{\mathup{cem}}$ and the reference solution $ u_k^{\mathup{fe}}$ are obtained. To evaluate the efficiency of the multiscale method, we consider the relative $ L^2 $ error and energy error between $ u_{k}^{\mathup{cem}}$ and $u_{k}^{\mathup{fe}} $ defined as
\[ E_L^k\coloneqq \frac{\norm{{u^{\mathup{fe}}_{k}} - u^{\mathup{cem}}_{k}}_{{L^2}(\Omega )}}{\norm{{u^{\mathup{fe}}_{k}}}_{{L^2}(\Omega )}}\ \text{and}\ E_a^k\coloneqq \frac{\norm{{u^{\mathup{fe}}_{k}} - u^{\mathup{cem}}_{k}}_a}{\norm{{u^{\mathup{fe}}_{k}}}_a},\]
and the iterative rate of $ u_{k-1}^{\mathup{cem}}$ and $u_{k}^{\mathup{cem}} $ in two norms defined as
\[ T_L^k\coloneqq \frac{\norm{{u^{\mathup{cem}}_{k-1}} - u^{\mathup{cem}}_{k}}_{{L^2}(\Omega )}}{\norm{{u^{\mathup{cem}}_{k-1}}}_{{L^2}(\Omega )}}\ \text{and}\ T_a^k\coloneqq \frac{\norm{{u^{\mathup{cem}}_{k-1}} - u^{\mathup{cem}}_{k}}_a}{\norm{{u^{\mathup{cem}}_{k-1}}}_a}.\]

\subsection{Model problem 1}
\begin{figure}
	\centering
	\begin{subfigure}[b]{0.24\textwidth}
		\centering
		\includegraphics[width=\textwidth]{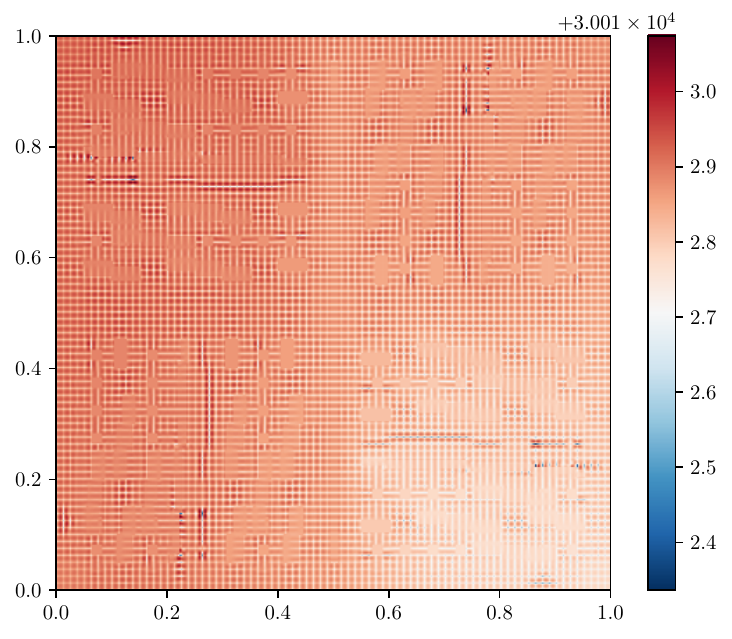}
		\caption{}\label{fig:C_Af1u0_iter0_imshow}
	\end{subfigure}
	\hfill
	\begin{subfigure}[b]{0.24\textwidth}
		\centering
		\includegraphics[width=\textwidth]{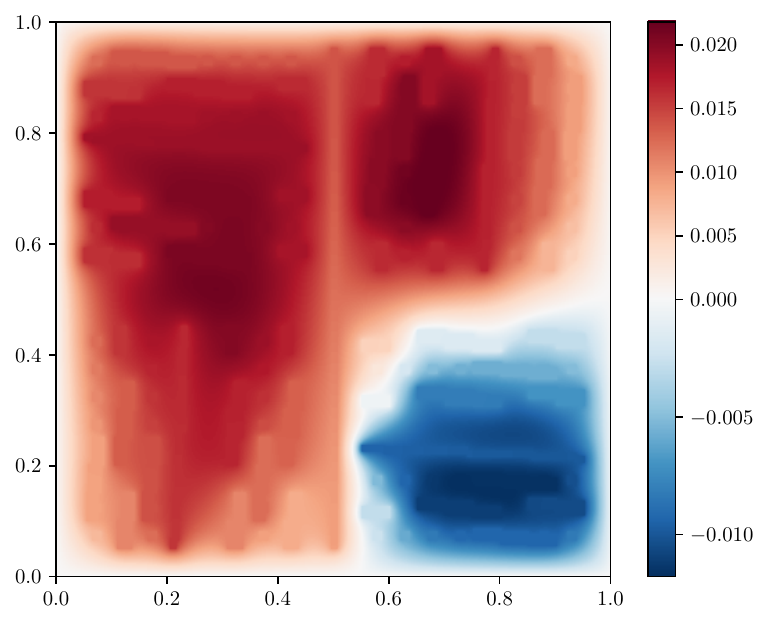}
		\caption{}
		\label{fig:C_Af1u0_iter1_imshow}
	\end{subfigure}
	\hfill
	\begin{subfigure}[b]{0.24\textwidth}
		\centering
		\includegraphics[width=\textwidth]{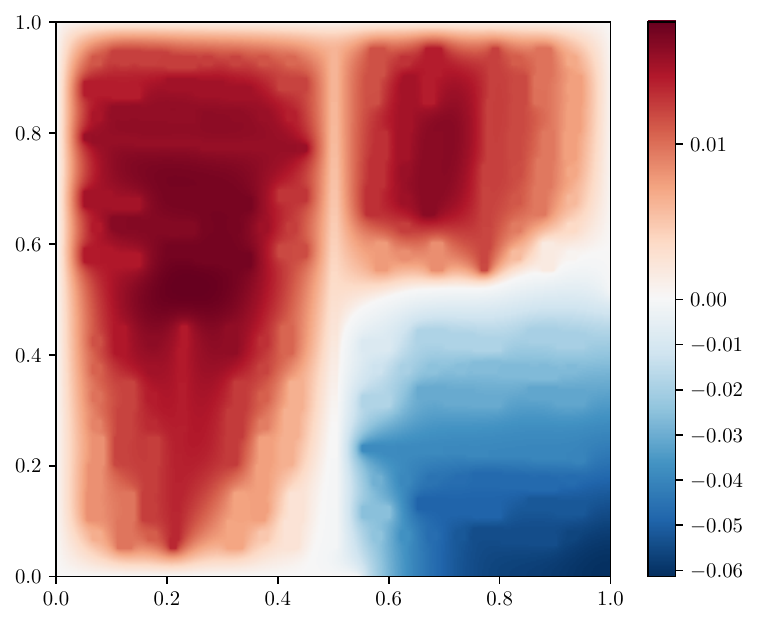}
		\caption{}
		\label{fig:C_Af1u0_iter2_imshow}
	\end{subfigure}
	\hfill
	\begin{subfigure}[b]{0.24\textwidth}
		\centering
		\includegraphics[width=\textwidth]{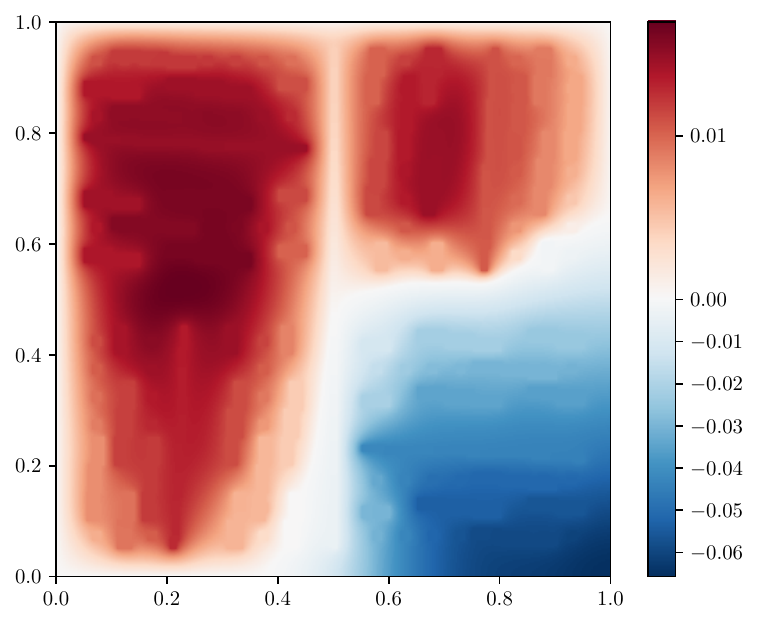}
		\caption{}
		\label{fig:C_Af1u0_iter7_imshow}
	\end{subfigure}
	\hfill
	\begin{subfigure}[b]{0.24\textwidth}
		\centering
		\includegraphics[width=\textwidth]{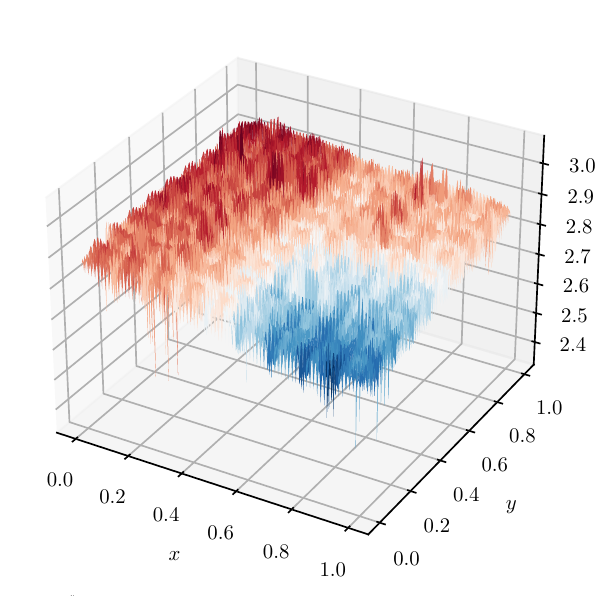}
		\caption{}\label{fig:C_Af1u0_iter0_3D}
	\end{subfigure}
	\hfill
	\begin{subfigure}[b]{0.24\textwidth}
		\centering
		\includegraphics[width=\textwidth]{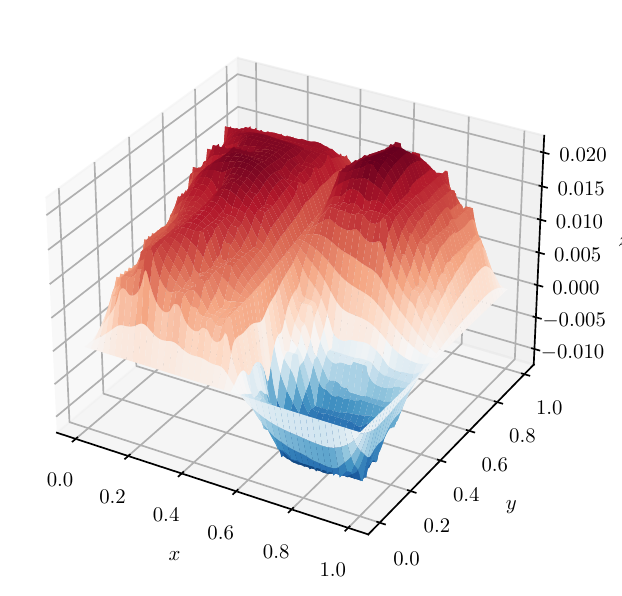}
		\caption{}
		\label{fig:C_Af1u0_iter1_3D}
	\end{subfigure}
	\hfill
	\begin{subfigure}[b]{0.24\textwidth}
		\centering
		\includegraphics[width=\textwidth]{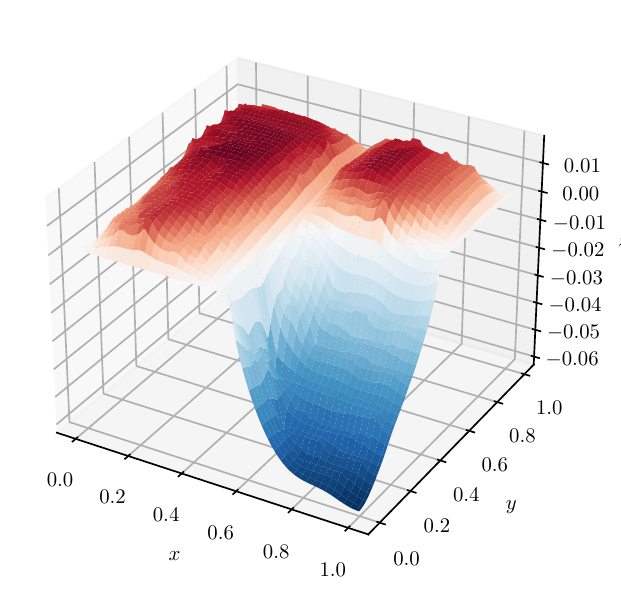}
		\caption{}
		\label{fig:C_Af1u0_iter2_3D}
	\end{subfigure}
	\hfill
	\begin{subfigure}[b]{0.24\textwidth}
		\centering
		\includegraphics[width=\textwidth]{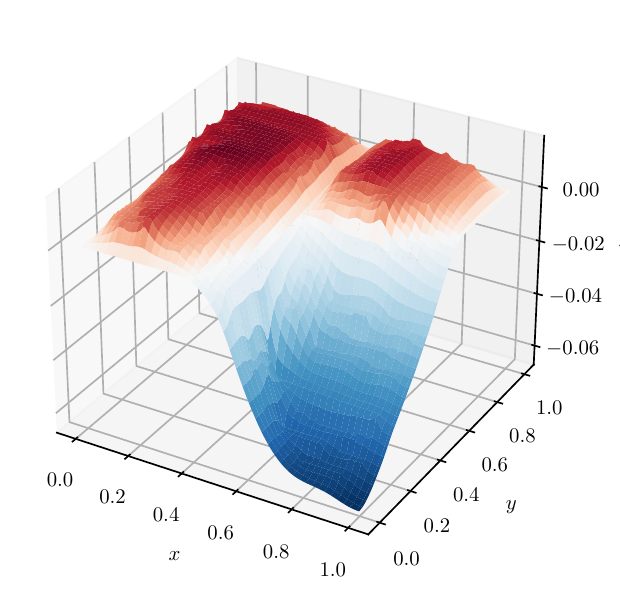}
		\caption{}
		\label{fig:C_Af1u0_iter7_3D}
	\end{subfigure}

	\caption{The solutions after iterations with medium A, using the source function $ f_1 $ and the initial value $u_0^\mathup{cem}=u_0^\mathup{fe}=u_0^0$. The first row shows the contour images, and the second row displays the 3D images. The iteration number: (a)(e)$ k $=1; (b)(f)$ k $=2; (c)(g)$ k $=3; (d)(h)$ k $=8.}
	\label{fig:C_Af1u0_results}
\end{figure}

\begin{figure}
	\centering
	\begin{subfigure}[b]{0.24\textwidth}
		\centering
		\includegraphics[width=\textwidth]{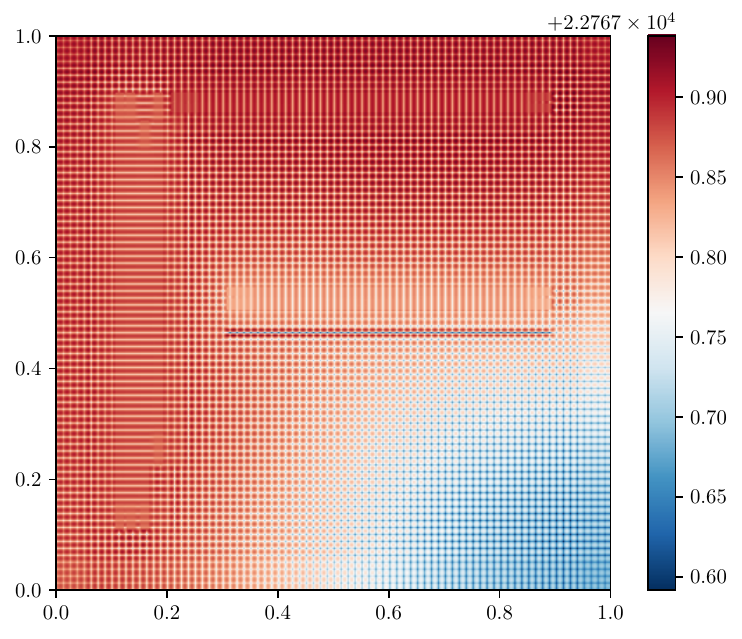}
		\caption{}\label{fig:C_Bf1u0_iter0_imshow}
	\end{subfigure}
	\hfill
	\begin{subfigure}[b]{0.24\textwidth}
		\centering
		\includegraphics[width=\textwidth]{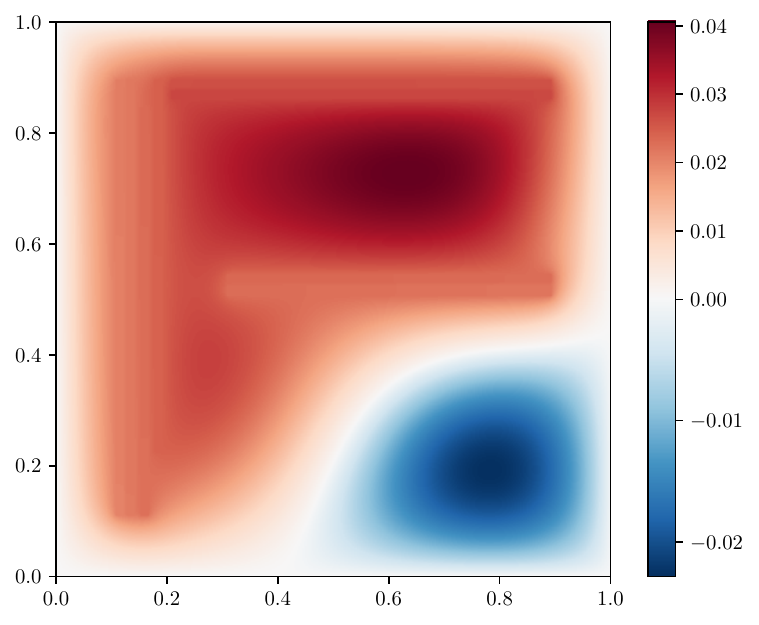}
		\caption{}
		\label{fig:C_Bf1u0_iter1_imshow}
	\end{subfigure}
	\hfill
	\begin{subfigure}[b]{0.24\textwidth}
		\centering
		\includegraphics[width=\textwidth]{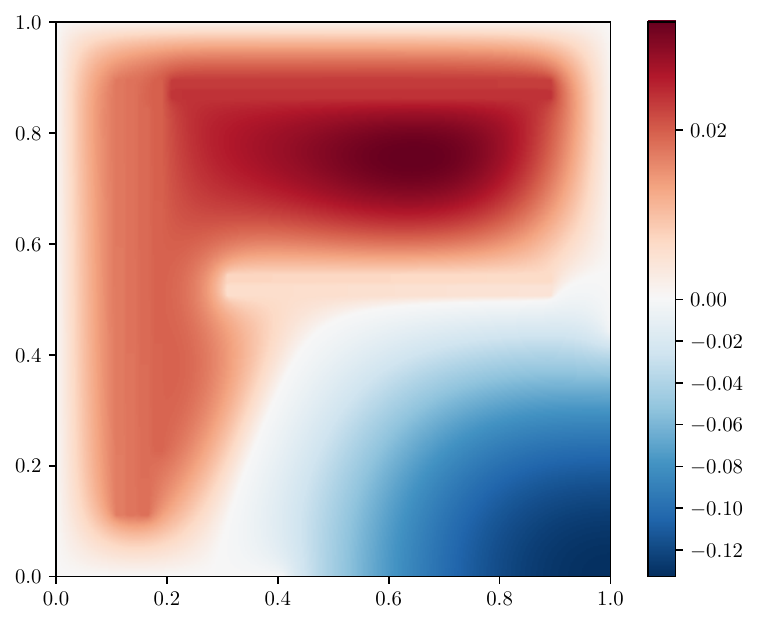}
		\caption{}
		\label{fig:C_Bf1u0_iter2_imshow}
	\end{subfigure}
	\hfill
	\begin{subfigure}[b]{0.24\textwidth}
		\centering
		\includegraphics[width=\textwidth]{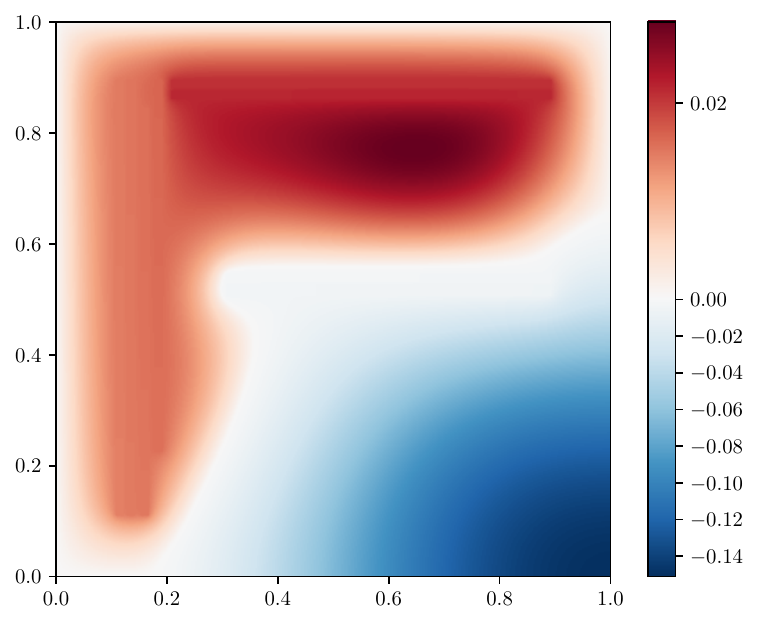}
		\caption{}
		\label{fig:C_Bf1u0_iter7_imshow}
	\end{subfigure}
	\hfill
	\begin{subfigure}[b]{0.24\textwidth}
		\centering
		\includegraphics[width=\textwidth]{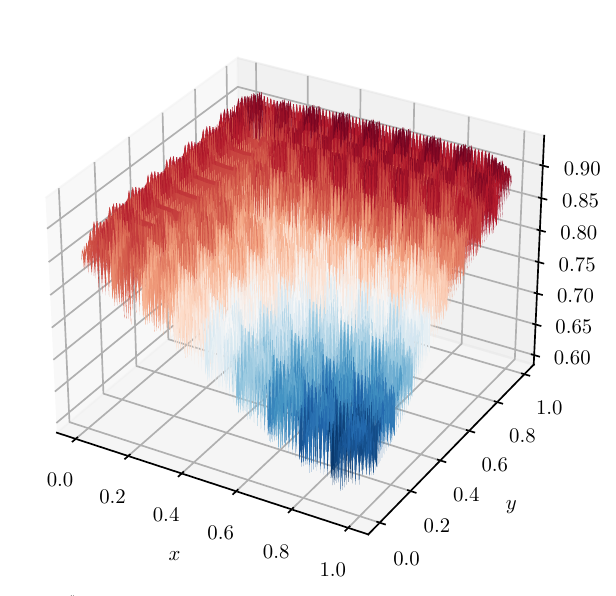}
		\caption{}\label{fig:C_Bf1u0_iter0_3D}
	\end{subfigure}
	\hfill
	\begin{subfigure}[b]{0.24\textwidth}
		\centering
		\includegraphics[width=\textwidth]{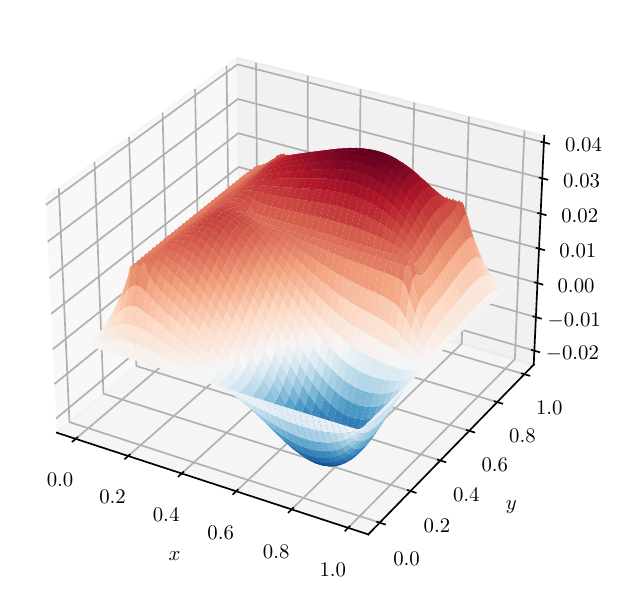}
		\caption{}
		\label{fig:C_Bf1u0_iter1_3D}
	\end{subfigure}
	\hfill
	\begin{subfigure}[b]{0.24\textwidth}
		\centering
		\includegraphics[width=\textwidth]{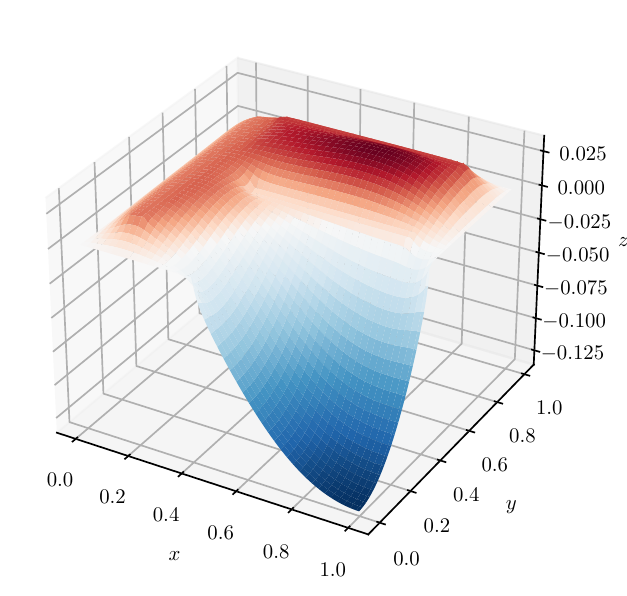}
		\caption{}
		\label{fig:C_Bf1u0_iter2_3D}
	\end{subfigure}
	\hfill
	\begin{subfigure}[b]{0.24\textwidth}
		\centering
		\includegraphics[width=\textwidth]{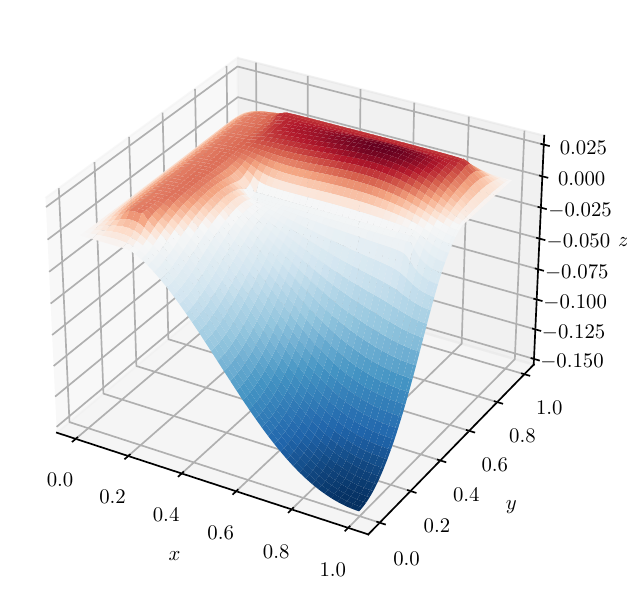}
		\caption{}
		\label{fig:C_Bf1u0_iter7_3D}
	\end{subfigure}

	\caption{The solutions after iterations with medium B, using the source function $ f_1 $ and the initial value $u_0^\mathup{cem}=u_0^\mathup{fe}=u_0^0$. The first row shows the contour images, and the second row displays the 3D images. The iteration number: (a)(e)$ k $=1; (b)(f)$ k $=2; (c)(g)$ k $=3; (d)(h)$ k $=8.}
	\label{fig:C_Bf1u0_results}
\end{figure}

We begin with a more challenging model by applying contact boundary conditions to all boundaries $ \partial \Omega $. While this might slightly different from the typical physics situation, our focus lies on the mathematical aspects and evaluating the effectiveness of our method with such nonlinear boundaries. The model problem is
\begin{equation}\label{eq:model1}
	\left\{
	\begin{aligned}
		 & - \Div \left( {\kappa \nabla u} \right) = f                                                    &  & \text { in } \Omega,          \\
		 & u \leqslant0, \quad \kappa\nabla u\cdot\bm{n} \leqslant0,\quad (\kappa\nabla u\cdot\bm{n})u =0 &  & \text { on } \partial \Omega, \\
	\end{aligned}
	\right.
\end{equation}
where $ \kappa$ is a highly heterogeneous permeability field generated from \cref{fig:Medium A} and \cref{fig:Medium B}, the source term $ f $ is taken as $ f_1 $ in the conducted experiments. Without loss of generality, the iteration processes are started with an initial value of $u_0^\mathup{cem}=u_0^\mathup{fe}=u_0^0=0$ in these two cases. To analyze the numerical results after iterations, we present contour images and 3D images respectively, see \cref{fig:C_Af1u0_results} and \cref{fig:C_Bf1u0_results}. With either of these two mediums, the results of the first numerical iteration oscillate, while after the second iteration, the values on the boundary converge to zero. Through a comparison of the outcomes from the iterations, it is apparent that the values of the contact boundaries have changed significantly, especially as the values decrease to less than or equal to zero. This result serves as an initial confirmation of the efficacy and accuracy of our methodology.

\subsection{Model problem 2}
In this subsection, we consider the following mixed contact boundary problem,
\begin{equation}\label{eq:model2}
	\left\{
	\begin{aligned}
		 & - \Div \left( {\kappa \nabla u} \right) = f                                                    &  & \text { in } \Omega, \\
		 & u=0                                                                                            &  & \text { on } \TD,    \\
		 & \kappa\nabla u\cdot\bm{n}=0                                                                    &  & \text { on } \TN,    \\
		 & u \leqslant0, \quad \kappa\nabla u\cdot\bm{n} \leqslant0,\quad (\kappa\nabla u\cdot\bm{n})u =0 &  & \text { on } \TC.    \\
	\end{aligned}
	\right.
\end{equation}
Let $\TD = [0, 1]\times{\{1\}},\  \TC =  [0, 1]\times{\{0\}}$, and $\TN = \partial\Omega \setminus (\TD\cup\TC )$.

\subsubsection{For testing accuracy}
To verify the effectiveness of our method for mixed contact boundary value problems, we consider highly heterogeneous permeability fields $\kappa$ generated from \cref{fig:Medium A} and \cref{fig:Medium B}. The source term for these cases is denoted as $f_2$. We initialize both iteration processes with $u_0^\mathup{cem} = u_0^\mathup{fe} = u_0^0$.

\begin{figure}
	\centering
	\begin{subfigure}[b]{0.24\textwidth}
		\centering
		\includegraphics[width=\textwidth]{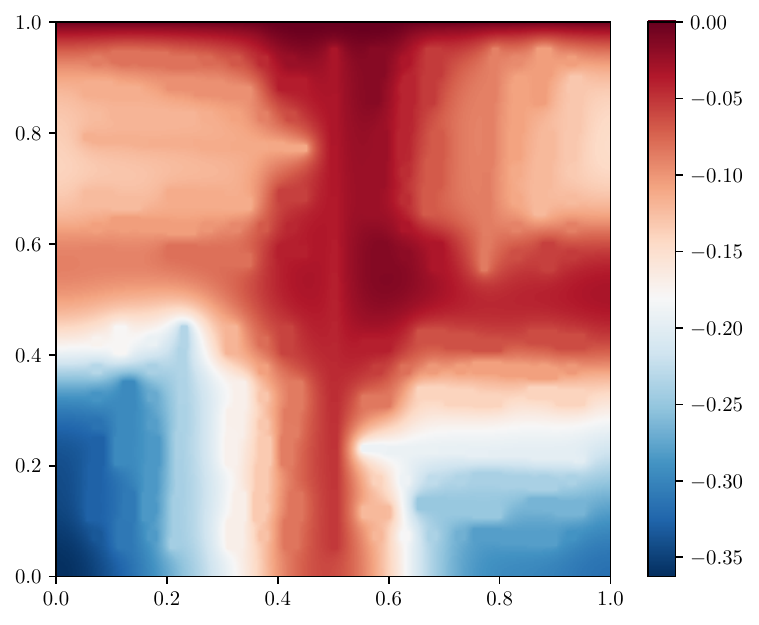}
		\caption{}\label{fig:dwC_Af2u0_iter0_imshow}
	\end{subfigure}
	\hfill
	\begin{subfigure}[b]{0.24\textwidth}
		\centering
		\includegraphics[width=\textwidth]{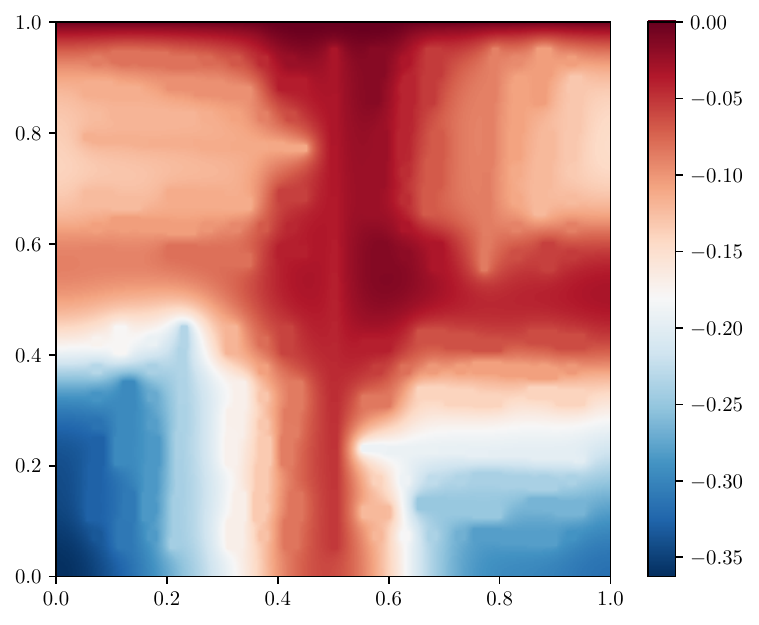}
		\caption{}
		\label{fig:dwC_Af2u0_iter1_imshow}
	\end{subfigure}
	\hfill
	\begin{subfigure}[b]{0.24\textwidth}
		\centering
		\includegraphics[width=\textwidth]{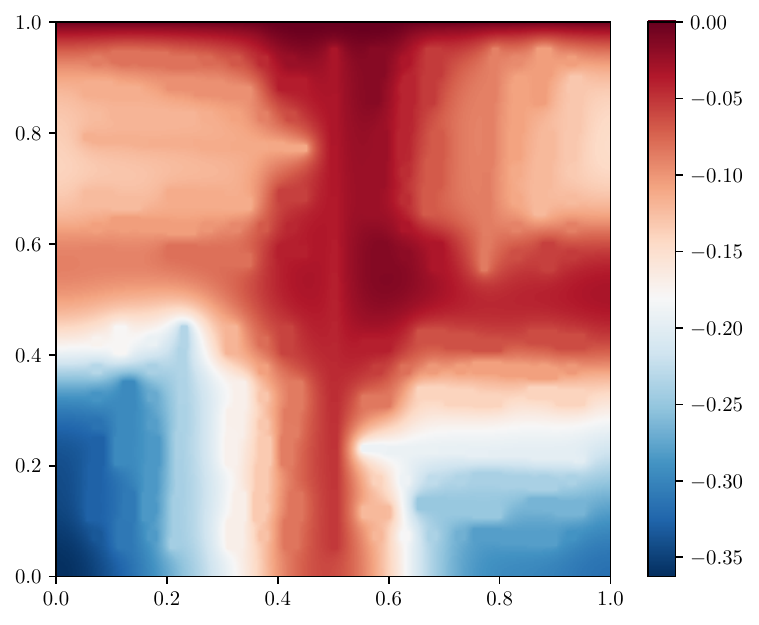}
		\caption{}
		\label{fig:dwC_Af2u0_iter2_imshow}
	\end{subfigure}
	\hfill
	\begin{subfigure}[b]{0.24\textwidth}
		\centering
		\includegraphics[width=\textwidth]{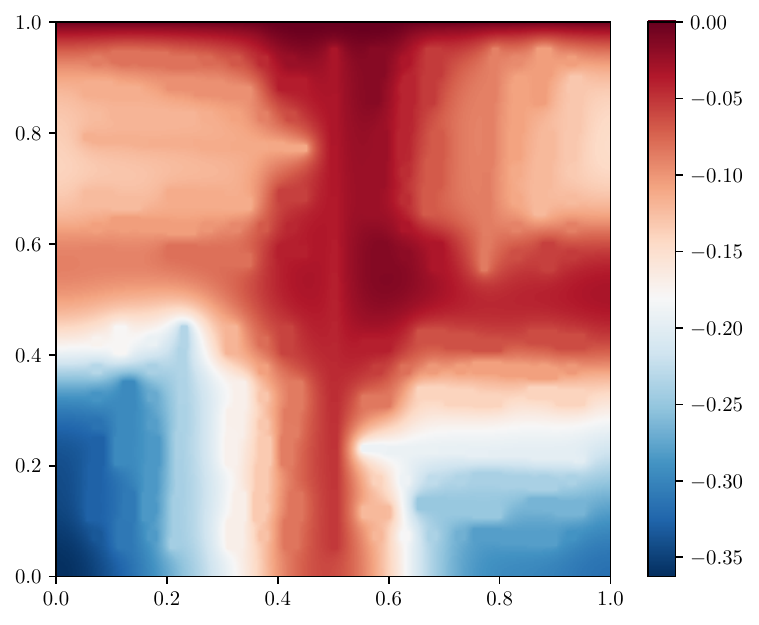}
		\caption{}
		\label{fig:dwC_Af2u0_iter7_imshow}
	\end{subfigure}
	\hfill
	\begin{subfigure}[b]{0.24\textwidth}
		\centering
		\includegraphics[width=\textwidth]{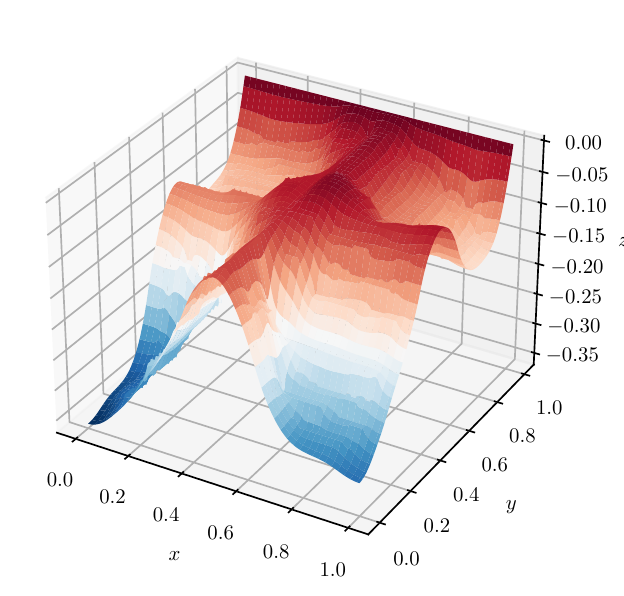}
		\caption{}\label{fig:dwC_Af2u0_iter0_3D}
	\end{subfigure}
	\hfill
	\begin{subfigure}[b]{0.24\textwidth}
		\centering
		\includegraphics[width=\textwidth]{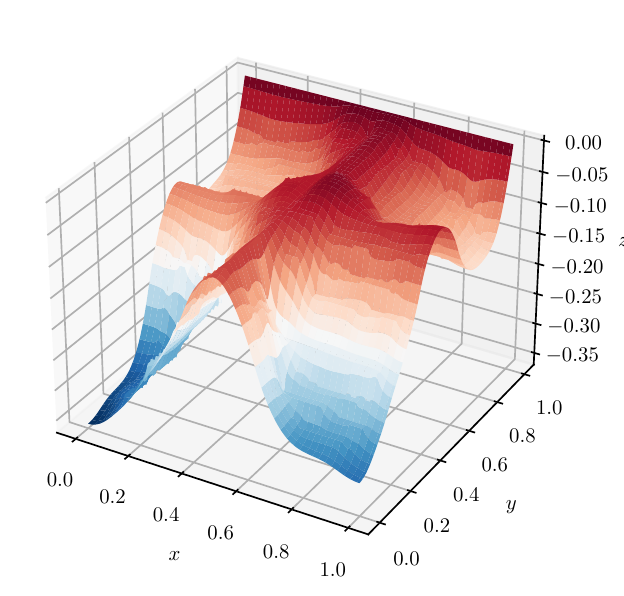}
		\caption{}
		\label{fig:dwC_Af2u0_iter1_3D}
	\end{subfigure}
	\hfill
	\begin{subfigure}[b]{0.24\textwidth}
		\centering
		\includegraphics[width=\textwidth]{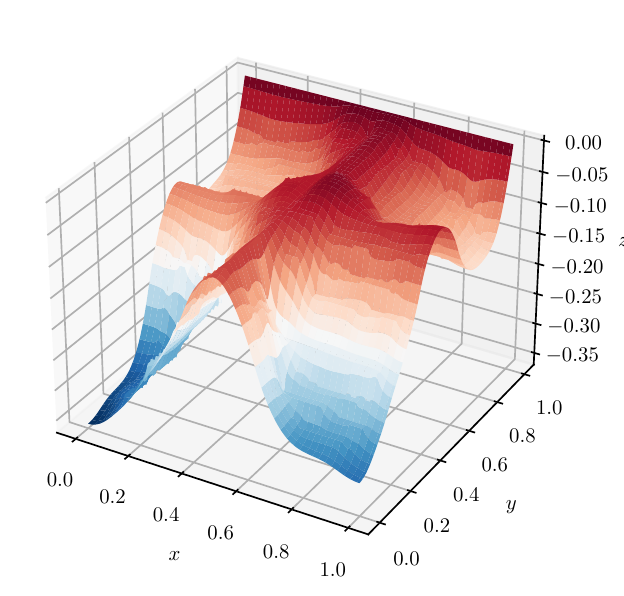}
		\caption{}
		\label{fig:dwC_Af2u0_iter2_3D}
	\end{subfigure}
	\hfill
	\begin{subfigure}[b]{0.24\textwidth}
		\centering
		\includegraphics[width=\textwidth]{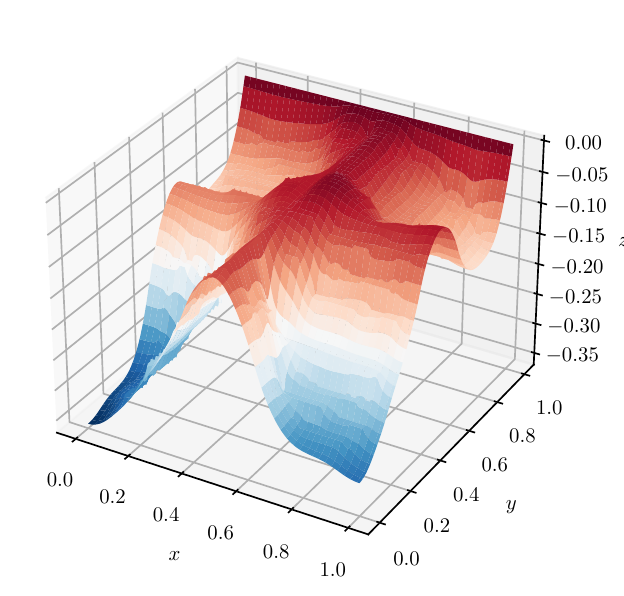}
		\caption{}
		\label{fig:dwC_Af2u0_iter7_3D}
	\end{subfigure}

	\caption{The solutions after iterations with medium A, using the source function $ f_2 $ and the initial value $u_0^\mathup{cem}=u_0^\mathup{fe}=u_0^0$. The first row shows the contour images, and the second row displays the 3D images. The iteration number: (a)(e)$ k $=1; (b)(f)$ k $=2; (c)(g)$ k $=3; (d)(h)$ k $=8.}
	\label{fig:dwC_Af2u0_results}
\end{figure}

\begin{figure}
	\centering
	\begin{subfigure}[b]{0.24\textwidth}
		\centering
		\includegraphics[width=\textwidth]{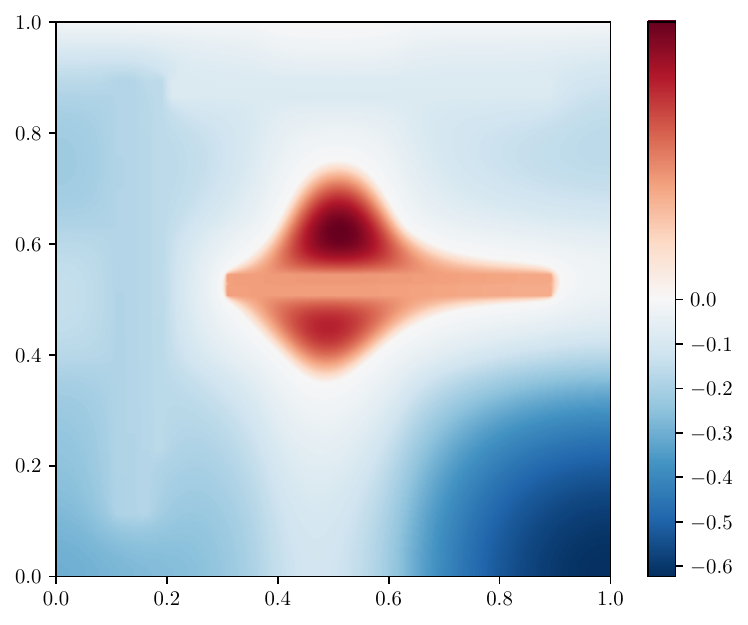}
		\caption{}\label{fig:dwC_Bf2u0_iter0_imshow}
	\end{subfigure}
	\hfill
	\begin{subfigure}[b]{0.24\textwidth}
		\centering
		\includegraphics[width=\textwidth]{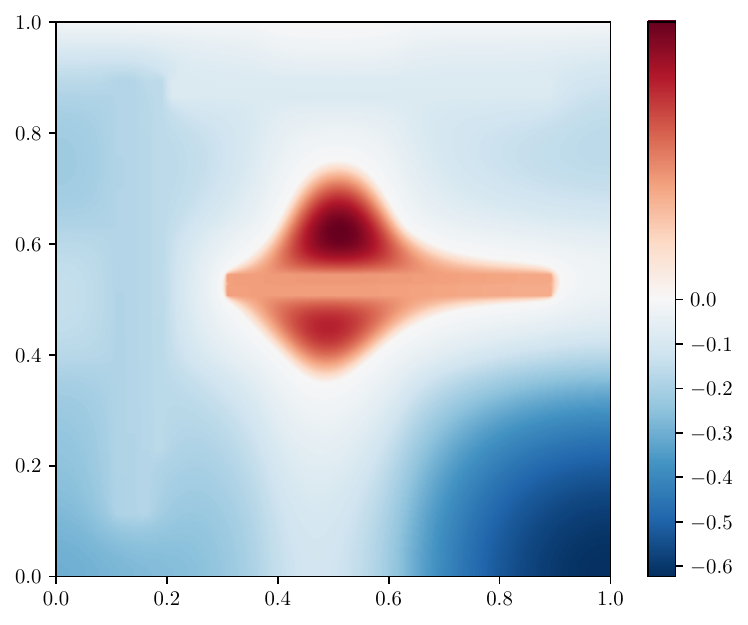}
		\caption{}
		\label{fig:dwC_Bf2u0_iter1_imshow}
	\end{subfigure}
	\hfill
	\begin{subfigure}[b]{0.24\textwidth}
		\centering
		\includegraphics[width=\textwidth]{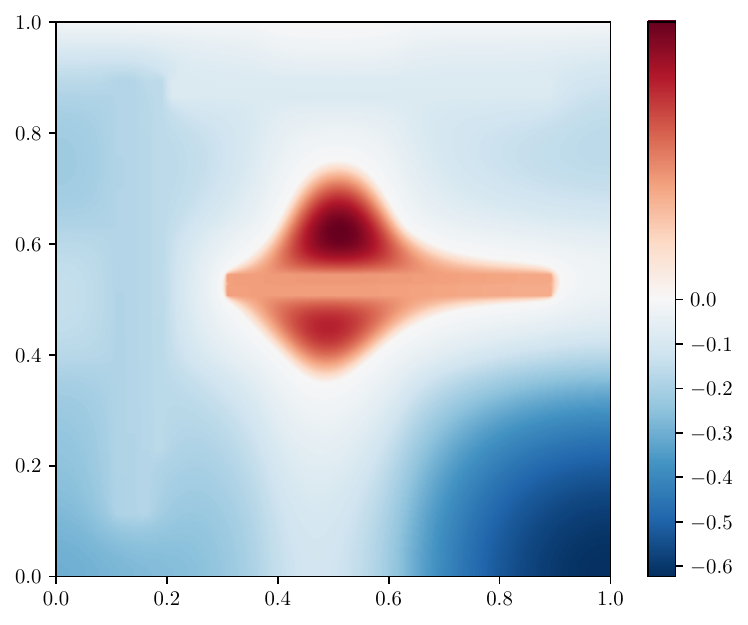}
		\caption{}
		\label{fig:dwC_Bf2u0_iter2_imshow}
	\end{subfigure}
	\hfill
	\begin{subfigure}[b]{0.24\textwidth}
		\centering
		\includegraphics[width=\textwidth]{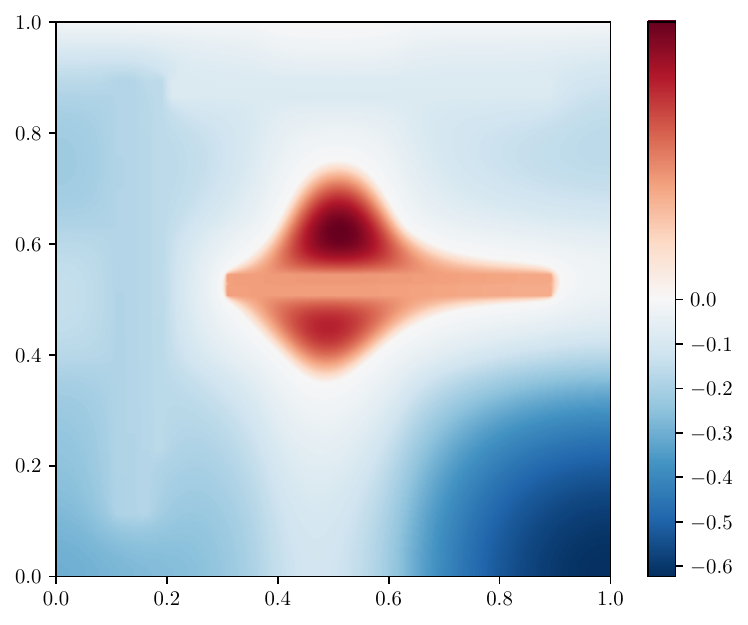}
		\caption{}
		\label{fig:dwC_Bf2u0_iter7_imshow}
	\end{subfigure}
	\hfill
	\begin{subfigure}[b]{0.24\textwidth}
		\centering
		\includegraphics[width=\textwidth]{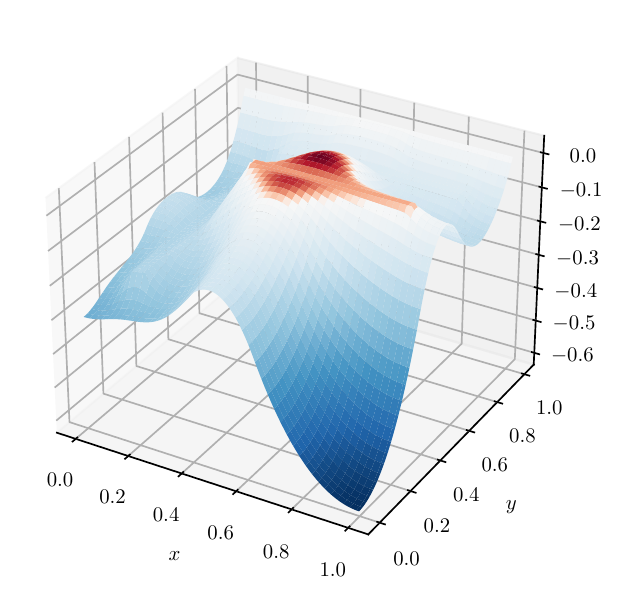}
		\caption{}\label{fig:dwC_Bf2u0_iter0_3D}
	\end{subfigure}
	\hfill
	\begin{subfigure}[b]{0.24\textwidth}
		\centering
		\includegraphics[width=\textwidth]{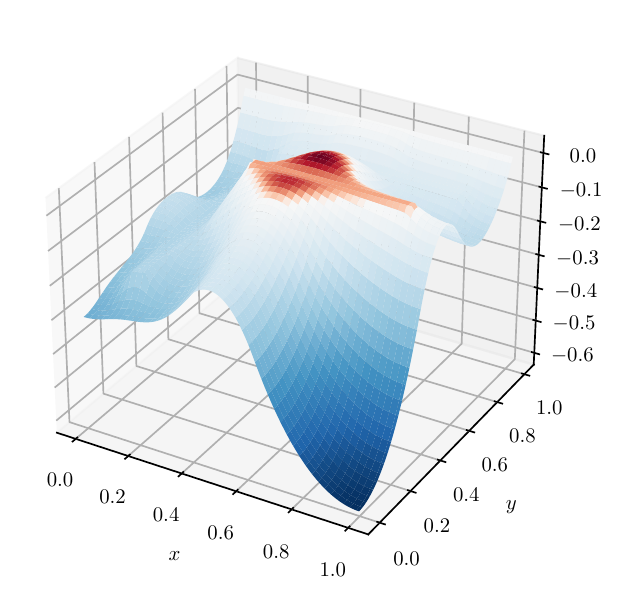}
		\caption{}
		\label{fig:dwC_Bf2u0_iter1_3D}
	\end{subfigure}
	\hfill
	\begin{subfigure}[b]{0.24\textwidth}
		\centering
		\includegraphics[width=\textwidth]{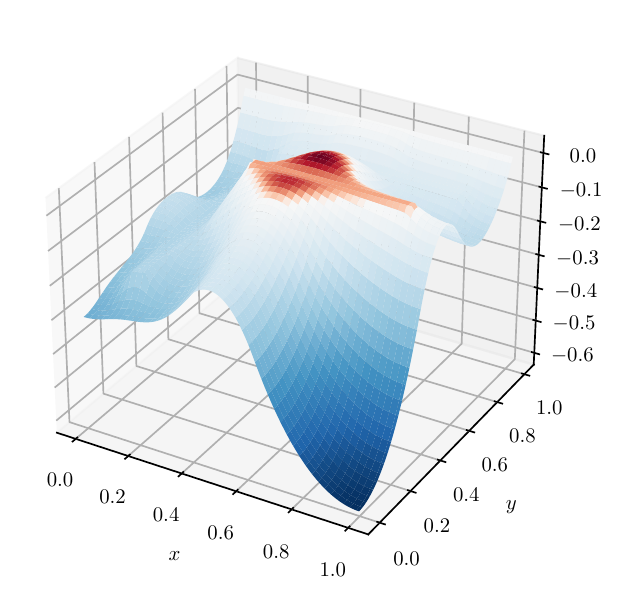}
		\caption{}
		\label{fig:dwC_Bf2u0_iter2_3D}
	\end{subfigure}
	\hfill
	\begin{subfigure}[b]{0.24\textwidth}
		\centering
		\includegraphics[width=\textwidth]{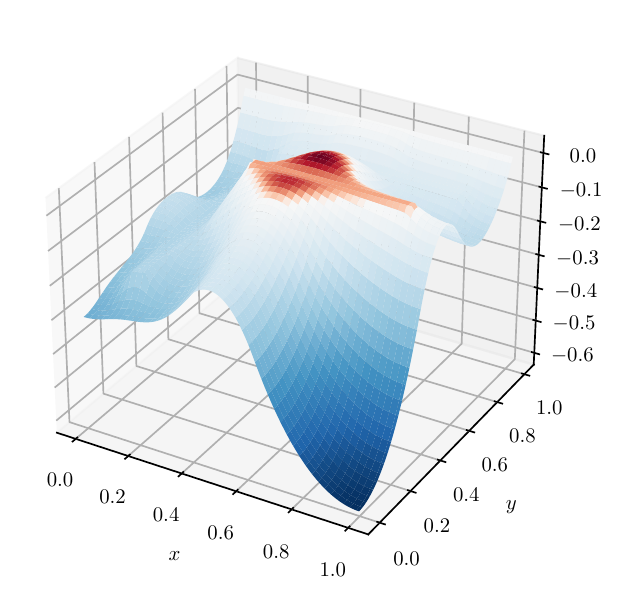}
		\caption{}
		\label{fig:dwC_Bf2u0_iter7_3D}
	\end{subfigure}

	\caption{The solutions after iterations with medium B, using the source function $ f_2 $ and the initial value $u_0^\mathup{cem}=u_0^\mathup{fe}=u_0^0$. The first row shows the contour images, and the second row displays the 3D images. The iteration number: (a)(e)$ k $=1; (b)(f)$ k $=2; (c)(g)$ k $=3; (d)(h)$ k $=8.}
	\label{fig:dwC_Bf2u0_results}
\end{figure}
The \cref{fig:dwC_Af2u0_results} and \cref{fig:dwC_Bf2u0_results} demonstrate that our proposed method effectively satisfies the conditions that the Dirichlet boundary values equal to 0 and the contact boundary values less than or equal to 0. Additionally, the results exhibit a rapid convergence rate when a suitable initial condition is utilized. These findings support the effectiveness of our method in handling the nonlinear and non-smooth conditions in the mixed contact boundary problem.

\subsubsection{For different initial guesses and parameters}\label{sec:initialandparameters}

To study the impact of various initial conditions, we set two more initial conditions with different power numbers on the results of numerical solution: $u_0^1=-x-y$ and $u_0^2=-x^2/2-y^2/2$. Also, we want to evaluate the impact of parameters on the results of solving nonlinear contact boundary problems. We select several sets of parameters and test the variation of the relative error and the iteration rate with serval iterations under these different parameter conditions. For these cases, the permeability field $\kappa$ is highly heterogeneous and is generated from \cref{fig:Medium C}. The source term is denoted as $f_3$. Through several experiments, we seek to analyze the behavior of different parameters under various initial conditions throughout the iterations.

\begin{figure}
	\centering
	\begin{subfigure}[b]{\textwidth}
		\centering
		\includegraphics[width=\textwidth]{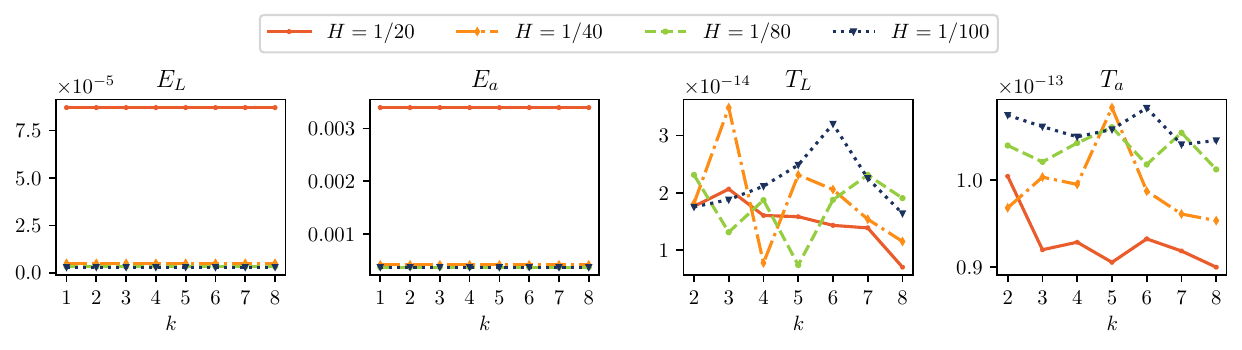}
		\caption{}
		\label{fig:Df3u00_1}
	\end{subfigure}
	\begin{subfigure}[b]{\textwidth}
		\centering
		\includegraphics[width=\textwidth]{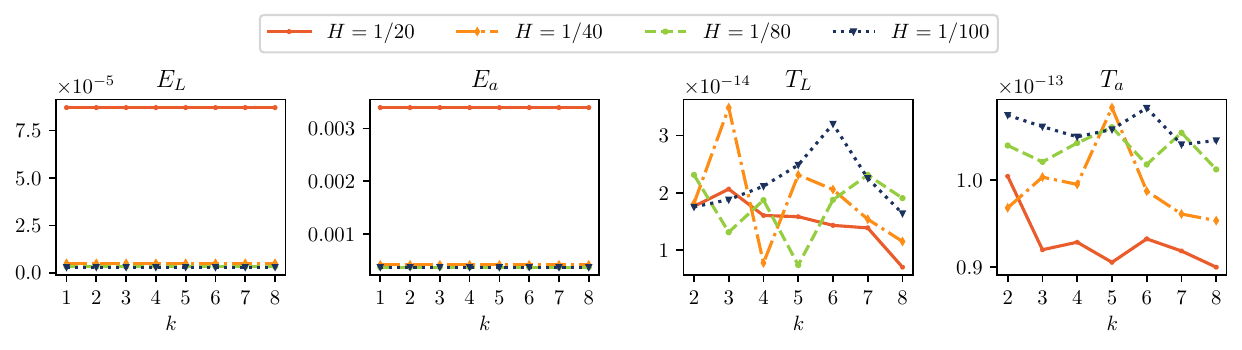}
		\caption{}
		\label{fig:Df3u01_1}
	\end{subfigure}
	\begin{subfigure}[b]{\textwidth}
		\centering
		\includegraphics[width=\textwidth]{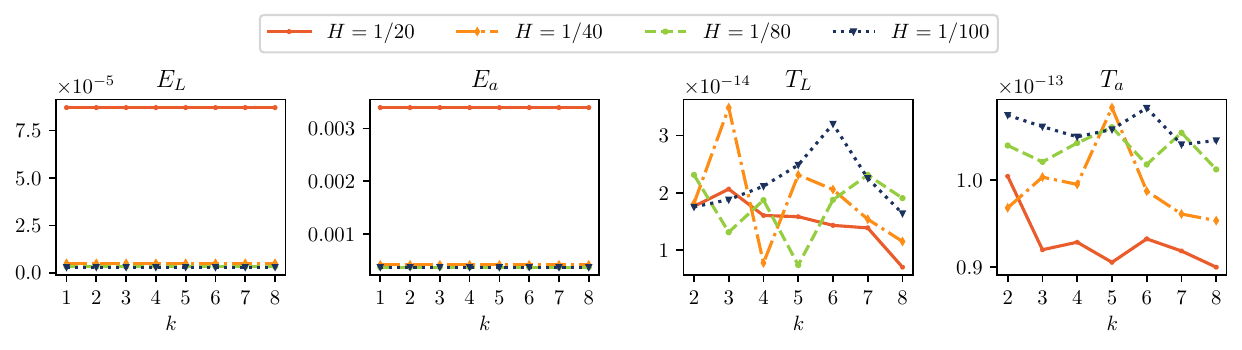}
		\caption{}
		\label{fig:Df3u02_1}
	\end{subfigure}
	\caption{The relative error and iterative rate of solutions under different coarse mesh sizes $ H $, with medium C, using the source function $ f_3 $. The initial conditions: (a) $ u_0^\mathup{cem}=u_0^\mathup{fe}=u_0^0 $, (b) $ u_0^\mathup{cem}=u_0^\mathup{fe}=u_0^1 $, (c) $u_0^\mathup{cem}=u_0^\mathup{fe}=u_0^2 $.}
	\label{fig:Df3initialu0_1}
\end{figure}

\begin{figure}
	\centering
	\begin{subfigure}[b]{\textwidth}
		\centering
		\includegraphics[width=\textwidth]{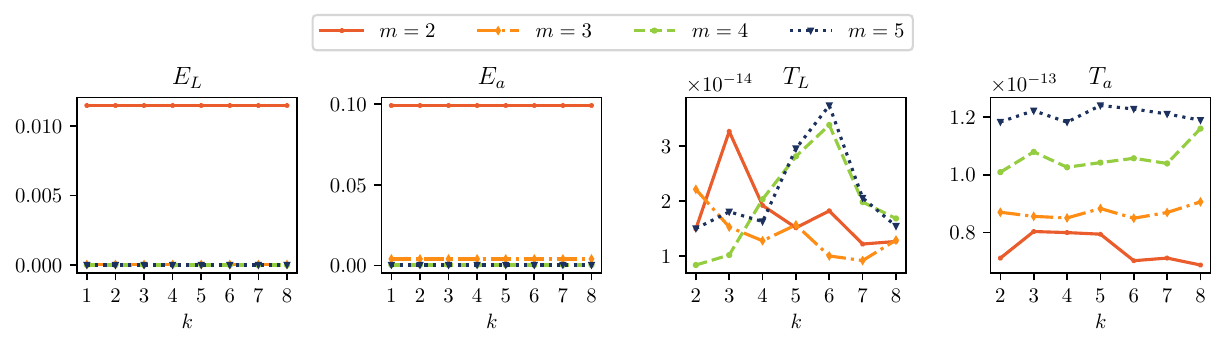}
		\caption{}
		\label{fig:Df3u00_5}
	\end{subfigure}
	\begin{subfigure}[b]{\textwidth}
		\centering
		\includegraphics[width=\textwidth]{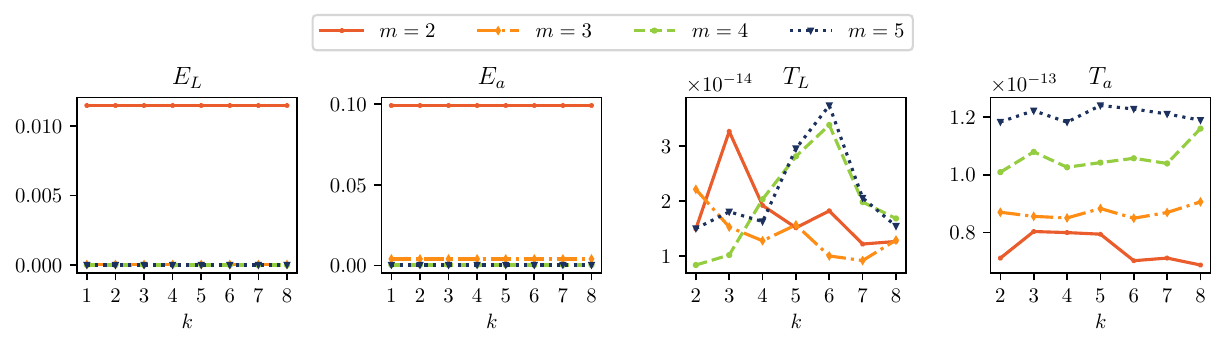}
		\caption{}
		\label{fig:Df3u01_5}
	\end{subfigure}
	\begin{subfigure}[b]{\textwidth}
		\centering
		\includegraphics[width=\textwidth]{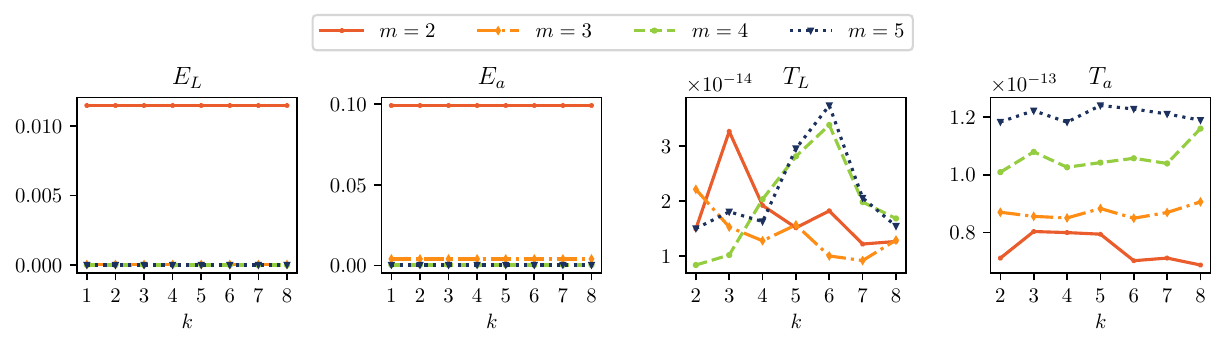}
		\caption{}
		\label{fig:Df3u02_5}
	\end{subfigure}
	\caption{The relative error and iterative rate of solutions under different oversampling layers $ m $, with medium C, using the source function $ f_3 $. The initial conditions: (a) $ u_0^\mathup{cem}=u_0^\mathup{fe}=u_0^0 $, (b) $ u_0^\mathup{cem}=u_0^\mathup{fe}=u_0^1 $, (c) $u_0^\mathup{cem}=u_0^\mathup{fe}=u_0^2 $.}
	\label{fig:Df3initialu0_5}
\end{figure}

\begin{figure}
	\centering
	\begin{subfigure}[b]{\textwidth}
		\centering
		\includegraphics[width=\textwidth]{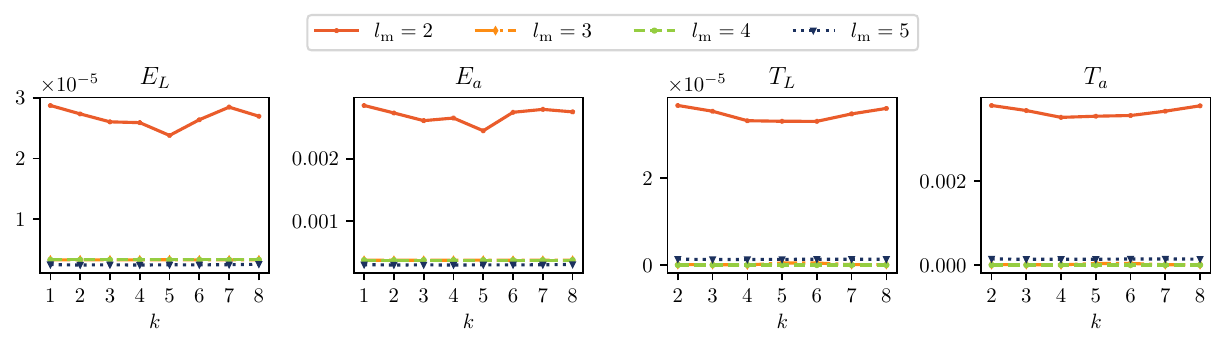}
		\caption{}
		\label{fig:Df3u00_4}
	\end{subfigure}
	\begin{subfigure}[b]{\textwidth}
		\centering
		\includegraphics[width=\textwidth]{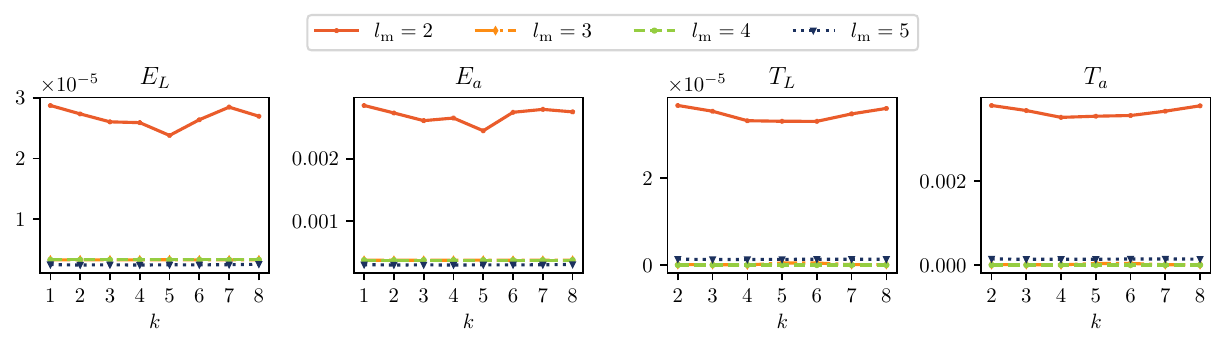}
		\caption{}
		\label{fig:Df3u01_4}
	\end{subfigure}
	\begin{subfigure}[b]{\textwidth}
		\centering
		\includegraphics[width=\textwidth]{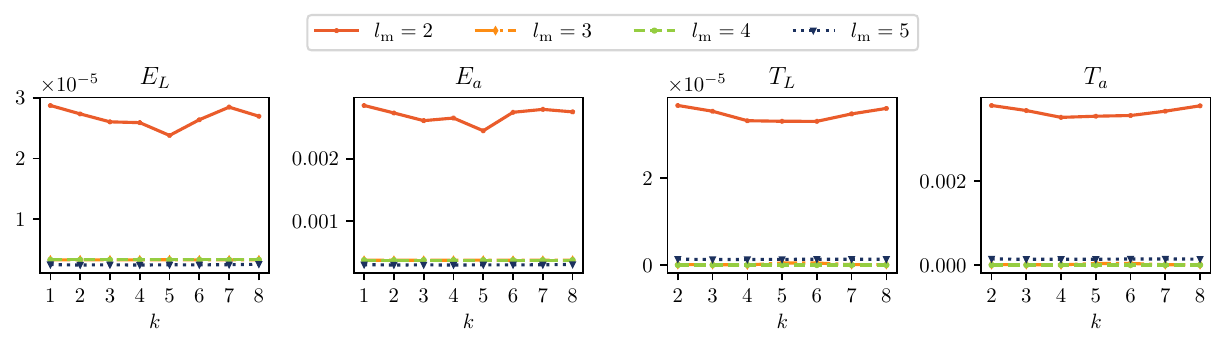}
		\caption{}
		\label{fig:Df3u02_4}
	\end{subfigure}
	\caption{The relative error and iterative rate of solutions under different eigenvector numbers $ l_\mathup{m}$, with medium C, using the source function $ f_3 $. The initial conditions: (a) $ u_0^\mathup{cem}=u_0^\mathup{fe}=u_0^0 $, (b) $ u_0^\mathup{cem}=u_0^\mathup{fe}=u_0^1 $, (c) $u_0^\mathup{cem}=u_0^\mathup{fe}=u_0^2 $.}
	\label{fig:Df3initialu0_4}
\end{figure}

The numerical experimental results, as depicted in \cref{fig:Df3initialu0_1}, \cref{fig:Df3initialu0_4} and \cref{fig:Df3initialu0_5}, provide valuable insights into our analysis. According to the definition of the iteration rate, only when $u_2^\mathup{cem}$ is solved, $T_L$ and $T_a$ can be calculated. Therefore, the plots start from $k=2$. In \cref{fig:Df3u00_1,fig:Df3u00_4,fig:Df3u00_5}, the initial condition $u_0^\mathup{cem}=u_0^\mathup{fe}=u_0^0$ is considered, and it shows the error and iterative rate. Similarly, \cref{fig:Df3u01_1,fig:Df3u01_4,fig:Df3u01_5}, as well as \cref{fig:Df3u02_1,fig:Df3u02_4,fig:Df3u02_5} illustrate the results with different initial conditions: $u_0^\mathup{cem}=u_0^\mathup{fe}=u_0^1$ and $u_0^\mathup{cem}=u_0^\mathup{fe}=u_0^2$, respectively. By comparing the three subplots in each figure, we observe that different initial conditions do not significantly affect the numerical solution obtained through the iterative process in this model setting.

The results indicate that as the number of iterations increases, the relative error remains consistently small, indicating that the difference between the solution obtained via our proposed method and the finite element solution is tiny. The relative error between the solutions obtained using our method and the Finite Element Method in the energy norm exhibits a similar trend to that in the $ L^2 $ norm, with the error in the energy norm consistently greater than that in the $L^2$ norm. 
The iterative rates keep being small. The consistent small iterative rates indicate that good convergence results are achieved after the first iteration in this model setting. The impact of further numerical iterations on the results is minimal, emphasizing our method is powerful. The numerical results demonstrate that our method achieves small errors and fast convergence rates across various sets of different parameters. However, different parameter settings still have some impact on the results.

In \cref{fig:Df3initialu0_1}, we investigate the influence of different coarse mesh sizes $H$, specifically $H=1/20$, $H=1/40$, $H=1/80$, and $H=1/100$. As expected, when $H=1/20$, the errors between the numerical solutions and the finite element solutions are the largest. As the coarse mesh size decreases, the error also decreases. Additionally, we observe that the rate of iteration in the numerical solution increases as the coarse mesh size decreases. By using a smaller coarse mesh size, the numerical solution can capture more detailed features of the problem, leading to reduced errors and improved iteration rates. However, it is important to make a balance between mesh size and computational resources, as using excessively fine meshes may lead to increased computational costs.

The \cref{fig:Df3initialu0_5} analyzes the effects of different oversampling layers $m$, varying from 2 to 5. We can observe that as the number of oversampling layers increases, the performance improves. The additional oversampling layers allow for finer resolution in capturing the details of the problem, leading to a more accurate representation of the solution. Actually, the relative error of ($ \kappa_\mathup{R} = 10^3$, $m = 3$) improves many times over ($ \kappa_\mathup{R} = 10^3$, $m = 2$). Similarly, as emphasized in the experiment of convergence rate, the $L^2$ norm errors ($\approx 10^{-14}$) are constantly smaller than energy norm errors ($\approx 10^{-13}$). This phenomenon reveals the potential of the CEM-GMsFEM in discovering contact models of high contrast problems.

Lastly, \cref{fig:Df3initialu0_4} explores the impact of different eigenvector numbers $l_\mathup{m}$, ranging from 2 to 5. Generally, increasing the eigenvector number leads to better results in the numerical solution. However, we also observe that when $l_\mathup{m}=3$ and $l_\mathup{m}=4$, the relative error and iteration rates are very similar. This suggests that the influence of the eigenvector number on the results may be affected by other parameters, such as coarse mesh sizes $H$ and oversampling layers $m$. Further analysis and experimentation are necessary to investigate the interactions between different parameters and their collective impact on the accuracy and convergence behavior of the numerical solution.

By examining these figures, we know how these parameters affect the accuracy and convergence of our method. This analysis allows us to select more suitable parameter settings for future experiments, thereby improving the overall performance of our approach.

\subsubsection{For interdependent parameters}
In this section, we examine how the numerical solution errors vary with two different parameters. The permeability fields $\kappa$ are highly heterogeneous and are generated based on the data provided in \cref{fig:Medium A} and \cref{fig:Medium B}. The source term is taken as $ f_1 $. Similar to previous experiments, we initiate the iteration process with $u_0^\mathup{cem}=u_0^\mathup{fe}=u_0^0$. Like previous experiment results in \cref{sec:initialandparameters}, the numerical results $ u^\mathup{cem}_{8} $ in this section remains stable to different parameters with very small error, therefore we consider it as the approximate solution $ u^\mathup{cem}_\infty $.

\begin{figure}
	\centering
	\begin{subfigure}[b]{0.8\textwidth}
		\centering
		\includegraphics[width=\textwidth]{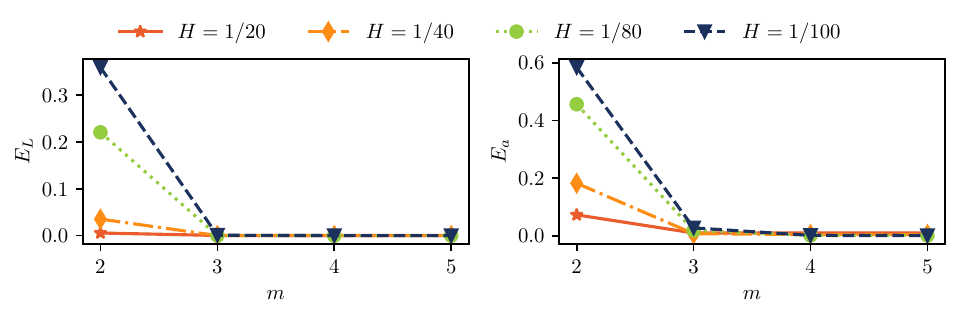}
		\caption{}
		\label{fig:AdwC_H_osly}
	\end{subfigure}
	\hfill
	\begin{subfigure}[b]{0.8\textwidth}
		\centering
		\includegraphics[width=\textwidth]{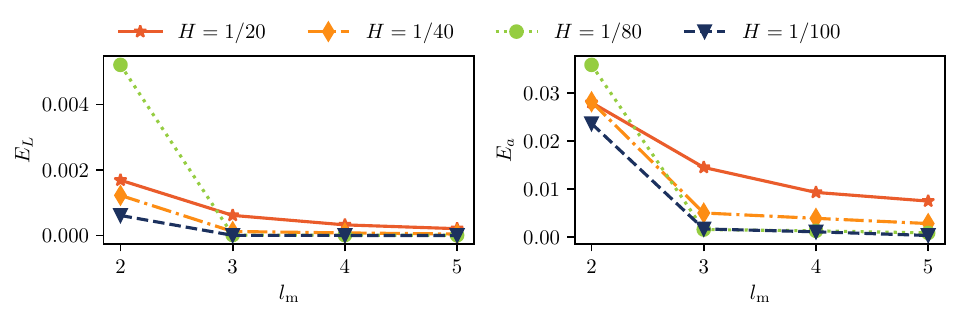}
		\caption{}
		\label{fig:AdwC_H_basis}
	\end{subfigure}
	\hfill
	\begin{subfigure}[b]{0.8\textwidth}
		\centering
		\includegraphics[width=\textwidth]{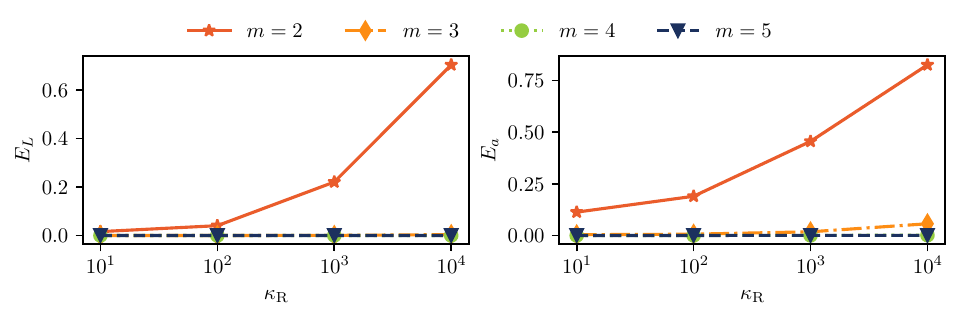}
		\caption{}
		\label{fig:AdwC_ctr_osly}
	\end{subfigure}
	\hfill
	\begin{subfigure}[b]{0.8\textwidth}
		\centering
		\includegraphics[width=\textwidth]{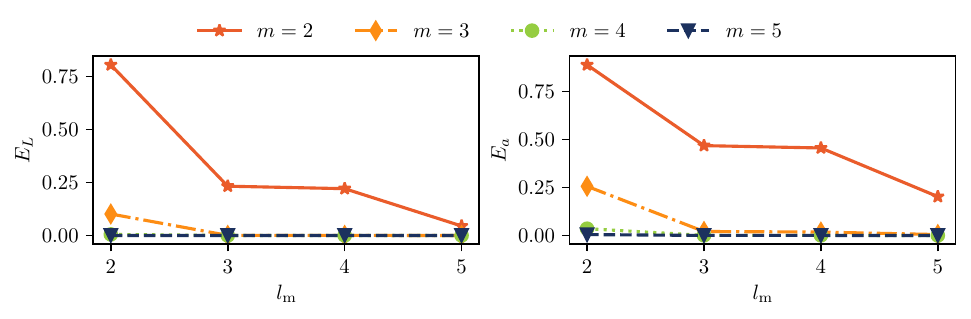}
		\caption{}
		\label{fig:AdwC_basis_osly}
	\end{subfigure}
	\caption{The relative errors in the energy and $ L^2 $ norm of numerical solutions with medium A, using the source function $f_1$ and the initial value $u_0^\mathup{cem}=u_0^\mathup{fe}=u_0^0$, (a)with different $ H $ and $ m $; (b)with different $ H $ and $ l_\mathup{m}$; (c)with different $ \kappa_\mathup{R}$ and $ m $; (d)with different $ l_\mathup{m}$ and $ m $.}
	\label{fig:Af1u0dwC_twopara}
\end{figure}

\begin{figure}
	\centering
	\begin{subfigure}[b]{0.8\textwidth}
		\centering
		\includegraphics[width=\textwidth]{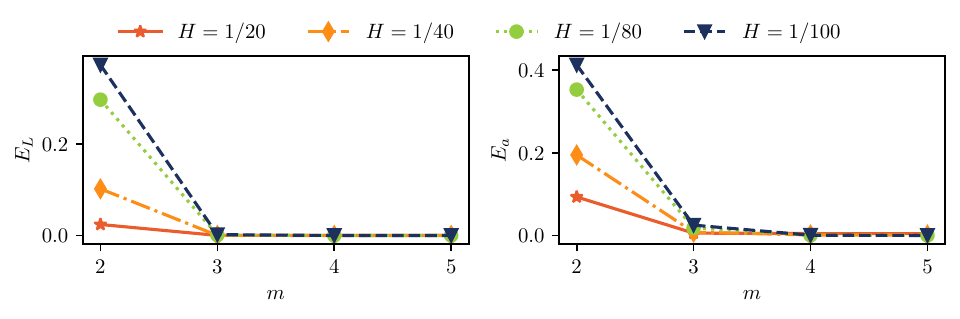}
		\caption{}
		\label{fig:BdwC_H_osly}
	\end{subfigure}
	\hfill
	\begin{subfigure}[b]{0.8\textwidth}
		\centering
		\includegraphics[width=\textwidth]{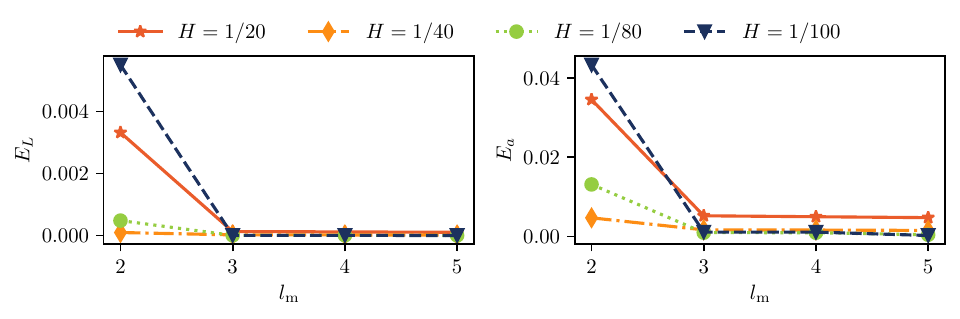}
		\caption{}
		\label{fig:BdwC_H_basis}
	\end{subfigure}
	\hfill
	\begin{subfigure}[b]{0.8\textwidth}
		\centering
		\includegraphics[width=\textwidth]{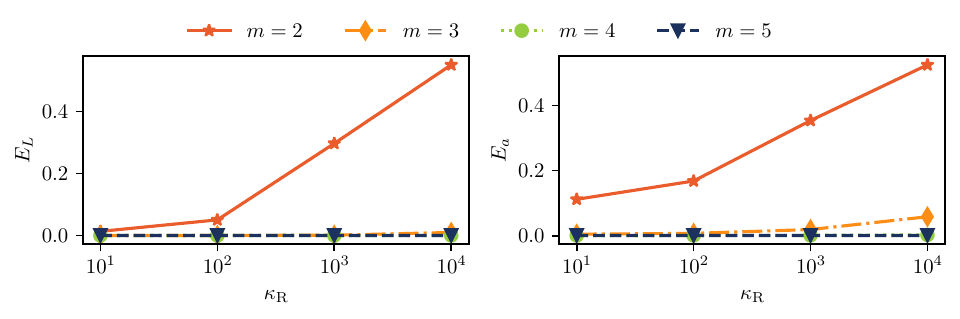}
		\caption{}
		\label{fig:BdwC_ctr_osly}
	\end{subfigure}
	\hfill
	\begin{subfigure}[b]{0.8\textwidth}
		\centering
		\includegraphics[width=\textwidth]{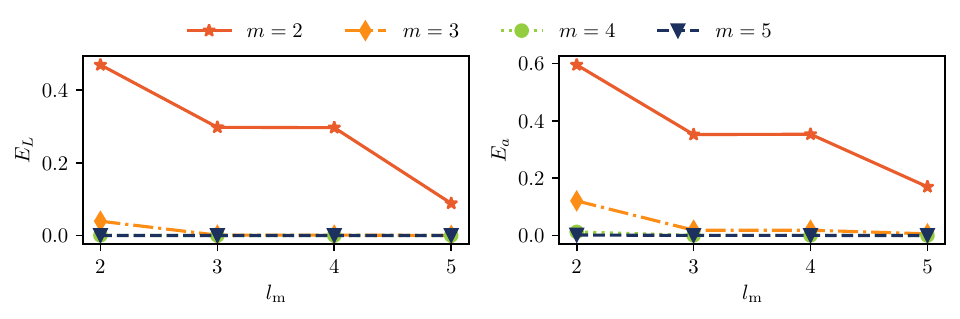}
		\caption{}
		\label{fig:BdwC_basis_osly}
	\end{subfigure}
	\caption{The relative errors in the energy and $ L^2 $ norm of numerical solutions with medium B, using the source function $f_1$ and the initial value $u_0^\mathup{cem}=u_0^\mathup{fe}=u_0^0$, (a)with different $ H $ and $ m $; (b)with different $ H $ and $ l_\mathup{m}$; (c)with different $ \kappa_\mathup{R}$ and $ m $; (d)with different $ l_\mathup{m}$ and $ m $.}
	\label{fig:Bf1u0dwC_twopara}
\end{figure}

The results are reported in \cref{fig:Af1u0dwC_twopara} and \cref{fig:Bf1u0dwC_twopara} for two different mediums. Notably, both sets of results sets of results have similar trends and characteristics. We will focus our analysis on the data of \cref{fig:Af1u0dwC_twopara}. We concentrate on investigating numerical errors by varying the coarse grids and oversampling layers, and then plot the results in \cref{fig:AdwC_H_osly}. An important observation is that errors will increase if we only reduce $ H $ while not enlarging oversampling layers $ m $, which is distinct from traditional finite element methods. Then, focusing on $ m = 2 $, we can see that $ E_L $ and $ E_a $ grows like $\bigO\left( {{H^{ - 1}}} \right) $. Another observation is that when $ m=3 $ and $ m=4 $, the error caused by different $ H $ becomes very tiny. A finer coarse mesh requires a larger oversampling layer to obtain a smaller relative error. In \cref{fig:AdwC_H_basis} illustrates how coarse grid and basis number affect the numerical error. The utilization of a finer mesh and an increased number of basis functions yields a substantial reduction in error. Notably, when the number of basis functions $ 3 $, the relative error can be effectively constrained to below $ 5\%$. We test numerical errors with different contrast ratios $ \kappa_\mathup{R}$ and oversampling layers $ m $, shown in \cref{fig:AdwC_ctr_osly}. It is not surprising that high contrast ratios will deteriorate numerical accuracy. Fortunately, the increase in oversampling layers $ m $ mitigates this deterioration. We then plot numerical error under different basis numbers and oversampling layers in \cref{fig:AdwC_basis_osly}. It is observed that the relative error diminishes as the number of basis functions increases and the oversampling domain expands. When increasing the number of oversampling layers to $ m=4 $ or $ m=5 $, the relative error remains small, consistently staying below 10\%. What' s more, further increasing the basis number $ l_m $ does not lead to significant improvements in the numerical results. Therefore, to strike a balance between computational complexity and result accuracy, we can choose the parameter values $ H=1/80 $, $ m=3 $, and $ l_m=3 $ for our experimental settings.

\section{Analysis}\label{sec:anal}
\begin{theorem}\label{Thm:Newton}
	Suppose that $ u_* $ is the real solution to the mixed contact boundary value problem \cref{eq:strong}.
	If $ \norm{u_0-u_*}_a $ is sufficiently small, then the semismooth Newton iteration is well-defined and converges superlinearly to $ u_* $.
\end{theorem}
\begin{proof}
	Since $ u_* $ is the real solution to the model problem \cref{eq:strong}, we have $ R(u_*)=0 $.	From \cref{subsec:semi-new}, we know that the residual functional $R(u)$ in \cref{eq:R(u)} is Newton differentiable on the open set $ U\subset V $, and the derivative of $R(u)$ is
	\[\left\langle {G(u)w,v} \right\rangle= a(w,v) + \frac{1}{\varepsilon }\int_{\TC }g'({u})wv\di\sigma.\]
	Since $ \int_{\TC }g'({u})w^2\di\sigma \geqslant 0$, we have
	\[\left\langle {G(u)w,w} \right\rangle \geqslant \norm{w}^2_a.\]
	Therefore, $ G(u) $ is nonsingular for all $ u \in U $ and $ \left\lbrace \norm{G(u)^{-1}}:u\in U \right\rbrace  $ is bounded.
	From Theorem 3.2 in \cite{Hintermueller2010}, the semismooth Newton method mentioned in \cref{subsec:semi-new} satisfies all conditions with the superlinear convergence to $ u_* $ if $ \norm{u_0-u_*}_a $ is sufficiently small.
\end{proof}

\begin{lemma}[See \cite{Ye2023}]\label{lemma:cem}
	Let $ u$ be the exact solution to the variational form \cref{eq:var_iter}, $ {\cal S}\{b(u),\kappa,f,p\} $ be the output numerical solution obtained in one iteration of CEM-GMsFEM solver with the input function $ u$. Then
	\begin{equation}\label{eq:lemma1}
		\norm{{\cal S}\{b(u),\kappa,f,p\}-u}_a\leqslant C'H\left( \norm{f}_{L^2} + \norm{p}_{L^2(\TN)} \right),
	\end{equation}
	where $ H $ is the coarse mesh size, $ C' $ is a constant that does not explicitly depend on the penalty parameter $ \varepsilon $ or the contrast ratio $ \kappa_\mathup{R} $ under carefully selecting the numbers of oversampling layers $ m $ and eigenvector number $ l_\mathup{m} $.
\end{lemma}

When considering two different input functions $ u' $ and $ u'' $, the main distinction in the outputs of the CEM-GMsFEM solver $ {\cal S}\{b(u'),\kappa,f,p\} $ and $ {\cal S}\{b(u''),\kappa,f,p\} $ are observed only on the contact boundary $ \TC $. Based on this, we may have the following assumption:
\begin{assumption}\label{assum}
	For the output of one step CEM-GMsFEM solver $ {\cal S}\{b,\kappa,f,p\} $, there exist positive constants $ r_0$ and $ C_* =1+\delta $ with $ \delta>0 $, such that for all input functions $ u'$ and $u''$ belong to a neighbourhood $\mathcal{B}_{r_0}(u_*) \subset V$, the following inequality holds:
	\[\norm{{\cal S}\{b(u'),\kappa,f,p\}-{\cal S}\{b(u''),\kappa,f,p\}}_a  \leqslant C_*\norm{u'-u''}_a.\]
\end{assumption}

\begin{theorem}\label{theorem:err}
	Under the \cref{assum}, suppose $ u_i\text{, }u_i^\mathup{cem} \in \mathcal{B}_{r_0}(u_*) \subset V$ for all $0\leqslant i\leqslant k $ and $ u_{0}=u_0^\mathup{cem}$. Let $ u_{i+1} $ be the exact solution derived by \cref{eq:var_iter}, $ u_{i+1}^\mathup{cem} $ represent the numerical solution derived by Iterative CEM-GMsFEM \cref{algo:cem} iterated from $u_{i}^\mathup{cem}$, $ C' $ and $ C_* $ be the constants defined in \cref{lemma:cem} and \cref{assum} respectively. Then there exist constants $ M $ and $ k_0 $ such that for all $ k\leqslant k_0 $,
	\begin{equation}\label{eq:thm_k+1_<1}
		\norm{u_{k+1}-u_{k+1}^\mathup{cem}}_a \leqslant C' MH\left( \norm{f}_{L^2} + \norm{p}_{L^2(\TN)} \right),
	\end{equation}
	where $ k_0=\left \lfloor{\frac{\log (M\delta +1)}{\log (\delta +1)}-1}\right \rfloor$.
	Furthermore, let $ r' $ be the convergence radius of semismooth Newton method and take $r\coloneqq \min\left\lbrace r_0\text{, }r'\right\rbrace$. The following estimate holds:
	\[
		\norm{u_*-u_{k+1}^\mathup{cem}}_a <2^{-(k+1)}r + C' MH\left( \norm{f}_{L^2} + \norm{p}_{L^2(\TN)} \right).\]
\end{theorem}
\begin{proof}
	Let $ \tilde{u}_{k+1}^{\mathup{cem}}$ represent the numerical solution obtained by Iterative CEM-GMsFEM \cref{algo:cem} iterated from the exact solution $ u_k $ to \cref{eq:var_iter}. By the definition, $ \tilde{u}_{k+1}^{\mathup{cem}} ={\cal S}\{b(u_k),\kappa,f,p\} $ and $ u_{k+1}^{\mathup{cem}} ={\cal S}\{b(u_{k}^{\mathup{cem}}),\kappa,f,p\}$. Because $ u_k$ and $u_k^\mathup{cem} $ belong to $ \mathcal{B}_{r_0}(u_*)$, we have
	\begin{equation}\label{eq:solvererror}
		\norm{\tilde{u}_{k+1}^{\mathup{cem}}-u_{k+1}^\mathup{cem}}_a  \leqslant C_*\norm{u_{k}-u_{k}^\mathup{cem}}_a
	\end{equation}
	by the \cref{assum}.
	We can split $\norm{u_{k+1}-u_{k+1}^\mathup{cem}}_a$ by introducing $ \tilde{u}_{k+1}^{\mathup{cem}} $:
	\[\begin{aligned}
			\norm{u_{k+1}-u_{k+1}^\mathup{cem}}_a
			 & =\norm{u_{k+1}-\tilde{u}_{k+1}^{\mathup{cem}}+\tilde{u}_{k+1}^{\mathup{cem}}-u_{k+1}^\mathup{cem}}_a,                     \\
			 & \leqslant \norm{u_{k+1}-\tilde{u}_{k+1}^{\mathup{cem}}}_a +\norm{\tilde{u}_{k+1}^{\mathup{cem}}-u_{k+1}^\mathup{cem}}_a . \\
		\end{aligned}\]
	For simplicity, we take $ e_k=u_{k}-u_{k}^\mathup{cem} $, which implies $ e_0=0 $.
	From \cref{lemma:cem} and \cref{assum}, it follows that
	\begin{equation}\label{eq:k+1cem}
		\norm{e_{k+1}}_a  \leqslant C' H\left( \norm{f}_{L^2} + \norm{p}_{L^2(\TN)} \right)+C_*\norm{e_k}_a.
	\end{equation}
	Dividing \cref{eq:k+1cem} on both sides by $ C_*^{k+1} $, we have
	\[\frac{\norm{e_{k+1}}_a}{C_*^{k+1}} \leqslant \frac{C' H}{C_*^{k+1}}\left( \norm{f}_{L^2} + \norm{p}_{L^2(\TN)} \right)+\frac{\norm{e_k}_a }{C_*^{k}}.\]
	Since $ u_i$ and $u_i^\mathup{cem} $ belong to $ \mathcal{B}_{r_0}(u_*) \subset V$ for all $0\leqslant i\leqslant k $, we have
	\[\frac{\norm{e_{i+1}}_a}{C_*^{i+1}}-\frac{\norm{e_i}_a }{C_*^{i}} \leqslant \frac{C' H}{C_*^{i+1}}\left( \norm{f}_{L^2} + \norm{p}_{L^2(\TN)} \right)\text{,  }\forall 0\leqslant i\leqslant k.\]
	We can sum the left-hand term for $0\leqslant i\leqslant k$ and build the inequality as follows:
	\[\begin{aligned}
			\frac{\norm{e_{k+1}}_a}{C_*^{k+1}}
			 & =\sum\limits_{i = 0}^k {\left( {\frac{\norm{e_{i+1}}_a }{C_*^{i+1}}-\frac{\norm{e_{i}}_a }{C_*^{i}}} \right)} \\
			 & \leqslant C' H\sum\limits_{i = 0}^k{\frac{1}{C_*^{i+1}}}\left( \norm{f}_{L^2} + \norm{p}_{L^2(\TN)} \right),  \\
			 & = C' H \frac{C_*^{k+1}-1}{C_*^{k+2}-C_*^{k+1}}\left( \norm{f}_{L^2} + \norm{p}_{L^2(\TN)} \right).            \\
		\end{aligned}\]
	As a result, we obtain
	\[	\begin{aligned}
			\norm{e_{k+1}}_a
			 & \leqslant C' H\frac{C_*^{k+1}-1}{C_*-1}\left( \norm{f}_{L^2} + \norm{p}_{L^2(\TN)} \right). \\
		\end{aligned}\]
	
	Since $ C_*=1+\delta $ with $ 0<\delta <1 $ and $ \delta $ small enough, then
	\[
		\norm{e_{k+1}}_a
		\leqslant C' H\frac{(1+\delta)^{k+1}-1}{(1+\delta)-1}\left( \norm{f}_{L^2} + \norm{p}_{L^2(\TN)} \right).\]
	Choosing a suitable constant $ M $ such that
	$$
		\frac{(1+\delta)^{k+1}-1}{(1+\delta)-1}\leqslant M,
	$$
	we have
	\[
		\norm{e_{k+1}}_a
		\leqslant C' MH\left( \norm{f}_{L^2} + \norm{p}_{L^2(\TN)} \right),\]
	for all $ k\leqslant k_0 $, where
	\[
		k_0=\left\lfloor {\frac{\log (\delta M+1)}{\log(1+\delta)}-1} \right\rfloor.\]

	From the \cref{Thm:Newton}, $ u_k $ converges superlinearly to $ u_* $, without loss of generality, assume $ \norm{u_*-u_{k}}_a  <2^{-k} r $ for all $ k\geqslant0 $.
	We can split $\norm{u_* - u_{k+1}^\mathup{cem}}_a$ and derive the following error estimates,
	\[\begin{aligned}
			\norm{u_*-u_{k+1}^\mathup{cem}}_a
			 & \leqslant \norm{u_*-u_{k+1}}_a +\norm{u_{k+1}-u_{k+1}^\mathup{cem}}_a,                      \\
			 & \leqslant  \norm{u_*-u_{k+1}}_a  + C'MH\left( \norm{f}_{L^2} + \norm{p}_{L^2(\TN)} \right), \\
			 & <2^{-(k+1)}r + C' MH\left( \norm{f}_{L^2} + \norm{p}_{L^2(\TN)} \right),                    \\
		\end{aligned}\]
\end{proof}

Indeed, if we choose appropriate initial $ u_0 $ and $ u_0^\mathup{cem} $, and make $ H $ sufficiently small, we can effectively control the error between the numerical solution and the exact solution at each iteration. We can also enhance the accuracy and convergence of the algorithm.

\section{Conclusion}
The paper presents a new iterative multiscale method based on the CEM-GMsFEM for solving the contact problem of Signorini type. This mixed contact boundary problem can be converted to a nonlinear unconstrained minimizing problem by applying the penalty method. The semismooth Newton method can be introduced to handle the challenges raised by the nonlinear and non-smooth contact conditions. Note that the variational form of an inhomogeneous Robin boundary value problem happen to be the same as that in semismooth Newton iteration. So the CEM-GMsFEM is employed iteratively to construct multiscale basis functions, which can effectively capture the solution behavior near the contact boundary. Numerical experiments were conducted using various heterogeneous coefficient profiles and tested with different initial conditions and parameters. The results validate our method holds particular significance in effectively addressing the challenges posed by non-linearity encountered in contact problems. The analysis demonstrates the proposed method is capable of providing accurate numerical solutions with fast convergence rates.

In future research, we plan to investigate the scalability and performance of large-scale problems by our method. Additionally, according to our experience, the construction of multiscale basis functions can be time-consuming. Therefore, it would be valuable to explore parallel computing approaches to speed up this process, which are rare in previous studies. Moreover, another interesting research direction is to extend our method to handle more realistic problems, such as dynamic contact and impact contact with friction.

\section*{Acknowledgment}

The research of Eric Chung is partially supported by the Hong Kong RGC General Research Fund (Project: 14305222).




\bibliographystyle{elsarticle-num}
\bibliography{cemgmsfem}
\end{document}